\pdfoutput=1
\RequirePackage{ifpdf}
\ifpdf 
\documentclass[pdftex]{sigma}
\else
\documentclass{sigma}
\fi

\newcommand\tensor{\otimes}
\newcommand\ST{\operatorname{ST}}

\newcommand\rz{{\color{red}0}}
\newcommand\ralpha{{\color{red}\alpha}}
\newcommand\rbeta{{\color{red}\beta}}
\newcommand\rgamma{{\color{red}\gamma}}
\newcommand\cS{\mathcal{S}}
\newcommand\bbP{\mathbb{P}}

\DeclareMathOperator{\pr}{pr}
\DeclareMathOperator{\dm}{dim}
\DeclareMathOperator{\gr}{gr}
\DeclareMathOperator{\id}{id}
\DeclareMathOperator{\img}{im}
\DeclareMathOperator{\homo}{End}
\DeclareMathOperator{\tf}{tf}
\DeclareMathOperator{\tr}{tr}

\DeclareMathOperator{\aut}{Aut}

\DeclareMathOperator{\sn}{span}

\DeclareMathOperator{\diag}{diag}
\DeclareMathOperator{\hm}{Hom}

\newcommand\fD{\mathfrak{D}}

\newcommand\cN{\mathcal{N}}
\newcommand\cL{\mathcal{L}}
\newcommand\GL{\operatorname{GL}}
\newcommand\Jac{\operatorname{Jac}}

\newcommand{\g}{\mathfrak{g}}
\newcommand{\gl}{\mathfrak{gl}}
\newcommand{\fh}{\mathfrak{h}}
\newcommand\bu{\mathbf{u}}
\newcommand\ad{\operatorname{ad}}
\newcommand\cG{\mathcal{G}}
\newcommand\Ad{\operatorname{Ad}}
\renewcommand\gr{\operatorname{gr}}

\newcommand\fa{\mathfrak{a}}
\newcommand\fann{\mathfrak{ann}}
\newcommand\ff{\mathfrak{f}}
\newcommand\fg{\mathfrak{g}}
\newcommand\finf{\mathfrak{inf}}
\newcommand\fM{\mathfrak{M}}
\newcommand\fp{\mathfrak{p}}
\newcommand\fS{\mathfrak{S}}
\newcommand\fs{\mathfrak{s}}
\newcommand\fsl{\mathfrak{sl}}
\newcommand\fU{\mathfrak{U}}

\newcommand\fm {\mathfrak{m}}
\newcommand\fz {\mathfrak{z}}
\newcommand\fX {\mathfrak{X}}
\newcommand\fq {\mathfrak{q}}

\newcommand{\R}{\mathbb{R}}
\newcommand{\C}{\mathbb{C}}
\newcommand{\Z}{\mathbb{Z}}
\newcommand{\V}{\mathbb{V}}

\newcommand\bbE{\mathbb{E}}
\newcommand\bbU{\mathbb{U}}

\newcommand\bbA{\mathbb{A}}
\newcommand\bbB{\mathbb{B}}

\newcommand\SL{\operatorname{SL}}

\newcommand\sfX{\mathsf{X}}
\newcommand\sfH{\mathsf{H}}
\newcommand\sfY{\mathsf{Y}}
\newcommand\sfZ{\mathsf{Z}}

\newcommand\cU{\mathcal{U}}
\newcommand\cA{\mathcal{A}}
\newcommand\cB{\mathcal{B}}
\newcommand\cE{\mathcal{E}}
\newcommand\cO{\mathcal{O}}
\newcommand\cW{\mathcal{W}}

\usepackage{tikz}
\usepackage{xfrac}
\usepackage{multirow}

\numberwithin{equation}{section}

\newtheorem{Theorem}{Theorem}[section]
\newtheorem{Corollary}[Theorem]{Corollary}
\newtheorem{Lemma}[Theorem]{Lemma}
\newtheorem{Proposition}[Theorem]{Proposition}
 { \theoremstyle{definition}
\newtheorem{Definition}[Theorem]{Definition}

\newtheorem{Remark}[Theorem]{Remark} }

\begin{document}
\allowdisplaybreaks

\renewcommand{\thefootnote}{}

\newcommand{\arXivNumber}{2301.09364}

\renewcommand{\PaperNumber}{058}

\FirstPageHeading

\ShortArticleName{On Uniqueness of Submaximally Symmetric Vector Ordinary Differential Equations}

\ArticleName{On Uniqueness of Submaximally Symmetric Vector\\ Ordinary Differential Equations of C-Class\footnote{This paper is a~contribution to the Special Issue on Symmetry, Invariants, and their Applications in honor of Peter J.~Olver. The~full collection is available at \href{https://www.emis.de/journals/SIGMA/Olver.html}{https://www.emis.de/journals/SIGMA/Olver.html}}}

\Author{Johnson Allen KESSY and Dennis THE}
\AuthorNameForHeading{J.A.~Kessy and D.~The}
\Address{Department of Mathematics and Statistics, UiT The Arctic University of Norway,\\ 9037 Troms\o, Norway}
\Email{\href{mailto:johnson.a.kessy@uit.no}{johnson.a.kessy@uit.no}, \href{mailto:dennis.the@uit.no}{dennis.the@uit.no}}

\ArticleDates{Received April 07, 2023, in final form August 01, 2023; Published online August 10, 2023}

\Abstract{The fundamental invariants for vector ODEs of order $\ge 3$ considered up to point transformations consist of generalized Wilczynski invariants and C-class invariants. An ODE of C-class is characterized by the vanishing of the former. For any fixed C-class invariant $\cU$, we give a local (point) classification for all submaximally symmetric ODEs of C-class with $\cU \not \equiv 0$ and all remaining C-class invariants vanishing identically. Our results yield generalizations of a well-known classical result for scalar ODEs due to Sophus Lie. Fundamental invariants correspond to the harmonic curvature of the associated Cartan geometry. A key new ingredient underlying our classification results is an advance concerning the harmonic theory associated with the structure of vector ODEs of C-class. Namely, for each irreducible C-class module, we provide an explicit identification of a lowest weight vector as a harmonic 2-cochain.}


\Classification{35B06; 53A55; 17B66; 57M60}

\renewcommand{\thefootnote}{\arabic{footnote}}
\setcounter{footnote}{0}

\section{Introduction}

Finite dimensionality of the contact symmetry algebra for scalar ODEs $u_{n+1}=f(t,u,u_1,\dots,u_n)$ of order $n+1\ge 4$
is a classical result due to Sophus Lie \cite{Lie1893} (see also \cite[Theorem~6.44]{Olver1995}). (We use jet notation $u_k$ instead of the more standard notation $u^{(k)}$ to denote the $k$-th derivative of~$u$ with respect to $t$.) The maximal symmetry dimension and the submaximal (i.e., next largest realizable) symmetry dimension are respectively
\begin{align*}
\fM:=n+5 \qquad \text{and}\qquad \fS := \begin{cases}
\fM-1 & \text{for}\ n=4\ \text{or}\ 6, \\
\fM-2 & \text{otherwise}.
\end{cases}
\end{align*}
 The former is realized locally uniquely by the trivial ODE $u_{n+1} = 0$. For ODEs realizing $\fS$, we have the following result (over $\C$) due to Lie \cite{Lie1924} (see also \cite[pp.~205--206]{Olver1995}): {\em Any submaximally symmetric scalar ODE of order $n+1 \ge 4$ is locally contact-equivalent to
\begin{itemize}\itemsep=0pt
\item[$(a)$] a linear equation, or
\item[$(b)$] exactly one of\,\footnote{In \cite[p.~206]{Olver1995}, the scalar ODE $3 u_2 u_4 - 5 (u_3)^2 = 0$ is also listed, but this is in fact contact-equivalent to $n u_{n-1} u_{n+1} - (n+1) (u_n)^2 = 0$ when $n=3$. We have verified this using Cartan-geometric techniques -- details will be given elsewhere.}
\begin{itemize}\itemsep=0pt
\item[$(i)$] $n=4\colon\ 9 (u_2)^2 u_5 - 45 u_2 u_3 u_4 + 40(u_3)^3 = 0$.
\item[$(ii)$] $n=6\colon\ 10(u_{3})^3u_{7}-70(u_{3})^2u_{4}u_{6} - 49(u_{3})^2(u_{5})^2 + 280u_{3}(u_{4})^2u_{5} -175(u_{4})^4 = 0$.
\item[$(iii)$] $n \neq 4,6\colon\ n u_{n-1} u_{n+1} - (n+1) (u_n)^2 = 0$.
\end{itemize}
\end{itemize}
}

The aim of our article is to establish analogous results for {\em vector} ODEs $\cE$ of order $n+1 \geq 3$:
\begin{align} \label{ODE}
\bu_{n+1} = \mathbf{f}(t,\bu,\bu_1,\dots,\bu_{n}),
\end{align}
where $\bu$ is an $\R^m$-valued function of $t$ (for $m \ge 2$), and $\bu_{k}$ is its $k$-th derivative. More precisely, we consider and completely resolve the classification problem (up to local contact equivalence) for submaximally symmetric vector ODEs \eqref{ODE} of order $\ge 3$ of {\em C-class} \cite{CDT2020, Cartan1938} (see below for motivation). Note that by the Lie--B\"{a}cklund theorem, contact-equivalence agrees with point-equivalence for vector ODEs.

For vector ODEs \eqref{ODE} of order $n+1 \ge 3$, the maximal and submaximal symmetry dimensions are
\begin{align}\label{E:MS}
\fM = m^2+(n+1)m+3 \qquad \text{and}\qquad \fS= \fM-2,
\end{align}
with the latter established in our earlier work \cite{KT2021}, along with numerous other symmetry gap results.
The trivial vector ODE $\bu_{n+1} = \mathbf{0}$ is locally uniquely maximally symmetric -- see for example \cite[Corollary~2.8]{KT2021}. Examples of some submaximally symmetric vector ODEs were given in \cite[Table 8]{KT2021}, but no definitive classification lists for the submaximal strata were asserted. This is a focus of our current article.

Following Cartan \cite{Cartan1938} (see also \cite{Bryant1991,CDT2020,Grossman2000}), a class of vector ODE \eqref{ODE} of order $\ge 3$ is said to be a {\em C-class} if it is invariant under all contact transformations, and all (contact) differential invariants of any ODE in this class are {\em first integrals} of that ODE. Hence, generic C-class equations (having sufficiently many functionally independent first integrals) can be solved using these invariants. In \cite[Theorem~4.1 and~4.2]{CDT2020}, the C-class was characterized by the vanishing of the {\em $($generalized$)$ Wilczynski invariants}. (This vanishing also leads to the existence of geometric structures on ODE solution spaces, which has been an important recent theme \cite{CDT2013, DT2006, GN2009, GN2010, Krynski2016}.) The Wilczynski invariants are a subset of the fundamental (relative) invariants (see Section~\ref{S: C-class}), which additionally consist of {\em C-class invariants} (in the terminology of \cite{KT2021}).

We note from \cite[Tables 8 and~10]{KT2021} that a vector ODE realizing $\fS$ given in \eqref{E:MS} is either a~3rd order ODE pair, i.e., $(n,m) = (2,2)$, of C-class or it is of {\em Wilczynski type} (i.e., an ODE with all C-class invariants vanishing identically). We will prove the following generalization of Lie's result above for vector ODEs:

\begin{Theorem} \label{T:I}
Any submaximally symmetric vector ODE \eqref{ODE} of order $n+1 \ge 3$ is either
\begin{itemize}\itemsep=0pt
\item[$(a)$] of Wilczynski type, or
\item[$(b)$] locally equivalent\footnote{More precisely, ``local equivalence'' here is meant in a neighbourhood of a point in $\cE$ where at least one of the C-class invariants is non-zero.} over $\R$ to exactly one of the three $3$rd order ODE pairs in Table {\rm \ref{Tab:CM}}. Over~$\C$, the two $3$rd order ODE pairs in the second row of Table {\rm \ref{Tab:CM}} are locally equivalent.
\end{itemize}
\end{Theorem}

Lie obtained his result for submaximally symmetric scalar ODEs using his complete classification of Lie algebras of contact vector fields on the (complex) plane and classified invariant ODEs having sufficiently many symmetries. Certainly, this approach generalizes to vector ODEs, but it is not feasible: complete classifications for Lie algebras of (point) vector fields on $\C^n$ or $\R^n$ for $n \ge 3$ are known to be very difficult to establish \cite{Doubrov2017,Schneider2018}. So, different techniques are needed to establish analogous results for submaximally symmetric vector ODEs.

Our approach to classifying all submaximally symmetric vector ODEs \eqref{ODE} of C-class of order $\ge 3$ is motivated by that of \cite{The2021, The2022} in the setting of parabolic geometries \cite{CS2009}, and is based on an equivalent reformulation of vector ODEs \eqref{ODE} as {\em $($strongly$)$ regular, normal Cartan ge\-ome\-tries}~$(\cG \to \cE, \omega)$ of type $(G,P)$ for a certain Lie group $G$ and closed subgroup $P \subset G$ \cite{CDT2020,Doubrov2001,DKM1999} (see Section~\ref{S:CG} below).

For such a (non-parabolic) Cartan geometry, the {\em harmonic curvature} $\kappa_H$, which corresponds to the fundamental invariants, is valued in a certain $P$-module that is completely reducible \cite[Corollary~3.8]{CDT2020}, so only the action of the reductive part $G_0 \subset P$ is relevant. Via a known algebraic Hodge theory associated with $G_0$, the codomain of $\kappa_H$ can be identified with a certain $G_0$-submodule $\bbE \subsetneq H^2(\g_-,\g)$ of a Lie algebra cohomology group called the {\em effective part} (see Definition \ref{D:E}). This has been already computed for ODEs \eqref{ODE} of order $3$ in \cite{Medvedev2010,Medvedev2011} and of order $\ge 4$ in \cite{DM2014}. The aforementioned fundamental invariants are valued in corresponding $G_0$-irreducible submodules $\bbU \subset \bbE$; see \cite[Table 6]{KT2021} for a summary. The irreducible C-class modules are listed in Table \ref{Tab:EP}.

We next formulate our second main result, which concerns the classification of vector ODEs \eqref{ODE} of C-class realizing the so-called constrained submaximal symmetry dimensions $\fS_\bbU$ identified in \cite[Table 2]{KT2021}. Fix an irreducible C-class module $\bbU = \bbB_4, \bbA_2^{\tr}, \bbA_2^{\tf} \subset \bbE$ (see Section~\ref{S: C-class}) and its corresponding C-class invariant $ \cU = \cB_4, \cA_2^{\tr}, \cA_2^{\tf}$ (see Section~\ref{S:CM}). Let $C_\cU$ denote the set of all ODEs \eqref{ODE} {\em with $\cU \not \equiv 0$ and all remaining C-class invariants vanishing identically} (equivalently, $0 \not \equiv \img(\kappa_H) \subset \bbU$), and let $\fS_\bbU$ denote the largest realizable symmetry dimension among ODEs in $C_\cU$. We will prove the following classification result:

\begin{Theorem} \label{T:Main1}
{\em Any} vector ODE \eqref{ODE} $\cE$ of C-class of order $n+1 \ge 3$ in $C_\cU$ realizing $\fS_\bbU$, near any point $x \in \cE$ with $\cU(x)\ne 0$, is {\em locally} $($point$)$ equivalent over $\R$ to exactly one of the ODEs given in Table {\rm \ref{Tab:CM}}. Over $\C$, the indicated $3$rd order ODEs for $\bbU = \bbB_4$ are locally equivalent.
\end{Theorem}

\begin{table}[h]\renewcommand{\arraystretch}{1.5}
\centering	
$\begin{array}{|c|c|c|c|} \hline
n& \begin{array}{c}
\text{Irreducible C-class}\\ \text{module } \bbU\subset \bbE
\end{array} &\fS_{\bbU} & \begin{array}{c}
\text{ODE of C-class with } 0\not\equiv \img(\kappa_H) \subset \bbU\\ \text{with symmetry dimension realizing } \fS_\bbU
\end{array} \\ \hline\hline
2 & \bbB_4 & \fM -m &\begin{array}{c}
\underset{(1 \le a \le m)}{u_3^a= \displaystyle\frac{3u_2^1 u_2^a}{2u^1_{1}}}
\end{array} \  \text{or}\  \begin{array}{c}
\underset{(1 \le a \le m)}{u_3^a= \displaystyle\frac{3u_1^1u_2^{1^{\vphantom{1}}} u_2^a}{1 + \bigl(u^1_1\bigr)^2}}\\
\end{array}\\\hline
\ge 3 & \bbA_2^{\tr}& \fM -m-1 & \begin{array}{c}
\underset{(1 \le a \le m)}{u_{n+1}^a= \displaystyle\frac{(n+1)u_n^{1^{\vphantom{1}}} u_n^a}{nu^1_{n-1}}}	\\
\end{array}\\	\hline		
\ge 2& \bbA_2^{\tf} & \fM -2m +1 +\delta_2^n &
\begin{array}{c}
\underset{(1 \le a \le m)}{u_{n+1}^a= \bigl(u_n^2\bigr)^2\delta_1^a}	\\
\end{array} \\ \hline
\end{array}$

\vspace{1mm}

(Recall $\fM = m^2+(n+1)m+3$ from \eqref{E:MS}.)

\caption{Classification over $\R$ of submaximally symmetric vector ODEs \eqref{ODE} of C-class of order $n+1 \geq 3$.}
\label{Tab:CM}
\end{table}

Our method for proving Theorems \ref{T:I} and \ref{T:Main1} will rely on the Cartan-geometric viewpoint for vector ODEs, and the associated computations will be efficiently done using representation theory. This will require important refinements to the existing structural results for vector ODEs of C-class stated in Table \ref{Tab:EP}. Such refinements constitute our final main result, which we now briefly describe. In our non-parabolic ODE setting, the aforementioned algebraic Hodge theory establishes a $G_0$-equivariant identification of $H^2(\g_-,\g)$ with the subspace ${\ker \square \subset \bigwedge^2 \g_-^* \tensor \g}$ of harmonic 2-cochains (see Section~\ref{S:CG}). Analogous to Kostant's theorem \cite{Kostant1961}, which is fundamental in the study of parabolic geometries, we may seek harmonic realizations of lowest weight vectors $\Phi_\bbU \in \bbU$ for each irreducible C-class submodule $\bbU \subset \bbE \subsetneq H^2(\g_-,\g)$. Our Theorem \ref{T:lwv} establishes such realizations (see Table \ref{Tab:lwv}). We anticipate that these structural results will be important for future geometric studies of the C-class and vector ODEs in general.

\section{Cartan geometries and vector ODEs of C-class} \label{S:CGP}

We briefly review the Cartan-geometric reformulation for vector ODEs \eqref{ODE} of order $\ge 3$ modulo point transformations, and summarize all relevant facts about vector ODEs of C-class.

\subsection{ODE geometry and symmetry}\label{S:GS}

We begin by summarizing \cite[Section~2.1]{KT2021}, which is based on \cite{Doubrov2001,DKM1999,DM2014}, and refer the reader to these articles for more details. The $(n+1)$-st order ODE \eqref{ODE} defines a submanifold ${\mathcal{E} = \{ \bu_{n+1} = \mathbf{f} \}}$ of co-dimension $m \ge 2$ in the space of $(n+1)$-jets of functions $J^{n+1}(\R, \R^m)$ that is transverse to the projection $\pi_n^{n+1}\colon J^{n+1}(\R, \R^m) \to J^{n}(\R, \R^m)$. Let $C$ denote the {\em Cartan distribution} on $J^{n+1}(\R, \R^m)$ with standard local coordinates $(t, \bu_0, \bu_1,\dots,\bu_{n+1})$, where $\bu_r = \big(u_r^1,\dots,u_r^m\big)$. Then $C$ is given by
\begin{align*}
C = \sn \{\partial_t + \bu_1 \partial_{\bu_0} + \dots + \bu_{n+1} \partial_{\bu_n},\, \partial_{\bu_{n+1}} \},
\end{align*}
where $\bu_i \partial_{\bu_j} :=\sum_{a=1}^m u_i^a \partial_{u_j^a}$ and $\partial_{\bu_r}$ refers to $\partial_{u_r^1},\dots,\partial_{u_r^m}$. We also consider the restriction of~$C$ to~$\cE$ and abuse notation by also referring to this distribution as $C$.
	
 Contact transformations are diffeomorphisms $\Phi\colon J^{n+1}(\R, \R^m) \to J^{n+1}(\R, \R^m)$ that preserve~$C$, i.e., ${\rm d}\Phi(C) = C$. By the Lie--B\"{a}cklund theorem, since $m \ge 2$, such transformations are the prolongations of diffeomorphisms on $J^0(\R, \R^m) \cong \R \times \R^m$, i.e., all such contact transformations are {\em point transformations}. Infinitesimally, a {\em point vector field} is a vector field $\xi \in \mathfrak{X}\bigl(J^{n+1}(\R,\R^m)\bigr)$ whose flow is a point transformation. Equivalently, $\cL_\xi C \subset C$, where $\cL_\xi$ is the Lie derivative with respect to $\xi$. A {\em point symmetry} of \eqref{ODE} is a point vector field that is tangent to $\cE$.

 We will consider ODEs \eqref{ODE} up to point transformations. The (point) geometry of $\cE$ is encoded by a pair $(E , V)$ of completely integrable sub-distributions of $C$ on $\cE$:
\begin{align} \label{E:EV}
E = \sn \left\{ \frac{\rm d}{{\rm d}t}:=\partial_t + \mathbf{u}_1\partial_{\mathbf{u}_0} + \cdots + \mathbf{u}_{n}\partial_{\mathbf{u}_{n-1}} + \mathbf{f}\partial_{\mathbf{u}_{n}}\right\}, \qquad
V= \sn \{\partial_{\mathbf{u}_n} \}.
\end{align}
(Note that integral curves of $E$ are lifts of solution curves to \eqref{ODE}.) Moreover, the distribution $D := E \oplus V \subset T\mathcal{E}$ is bracket-generating, and its weak-derived flag defines the following filtration on $T\cE$:
 \begin{align*}
 T\mathcal{E} = D^{-n-1} \supset \cdots \supset D^{-2} \supset D^{-1} := D,
 \end{align*}
 where $D^{-i-1} := D^{-i} + \big[D^{-i}, D^{-1}\big]$
for $i > 0$. Since $\big[\Gamma\bigl(D^j\bigr), \Gamma\bigl(D^k\bigr)\big] \subset \Gamma\bigl(D^{j+k}\bigr)$, then the pair $(\mathcal{E}, \{D^i\})$ forms a {\em filtered manifold}.
As we will describe below, this leads to the formulation of an ODE \eqref{ODE} as a filtered $G_0$-structure \cite[Section~2.1]{Andreas2017}.

Letting $T^i\cE := D^i \subset T\cE$ for $-n-1 \le i \le -1$ and $T^0\cE := 0$, we define
$\gr(T \mathcal{E}) := \bigoplus^{-1}_{i=-n-1} \gr_i(T\mathcal{E})$ where $\gr_i(T\mathcal{E}) := T^i\cE / T^{i+1}\cE$. Let $\gr_i(T_x\cE)$ denote the fiber of $\gr_i(T\cE)$ at $x \in \cE$, i.e., $\fm_i(x):=\gr_i(T_x\cE)= T^i_x\cE / T^{i+1}_x\cE$. Then $\fm (x):=\gr(T_x \cE )= \bigoplus^{-1}_{i=-n-1}\fm_i(x)$ is a~nilpotent graded Lie algebra (NGLA) under the (Levi) bracket induced by the Lie bracket of vector fields. It is called the {\em symbol algebra} at $x$. Since the symbol algebras at all points are isomorphic, then we let $\fm$ denote a fixed NGLA with $\mathfrak{m} \cong \fm (x)$, $\forall x \in \cE$, and we say that $(\cE, \{D^i\})$ is {\em regular of type} $\fm$.

Let $\aut_{\gr}(\fm) \le \GL(\fm)$ be the subgroup that preserves the grading of $\fm$. Since $\fm $ is generated by $\fm_{-1}$, then we have $\aut_{\gr}(\mathfrak{m}) \hookrightarrow\GL(\fm_{-1})$. For $x \in \cE$, we let $F_{\gr}(x)$ denote the set of all NGLA isomorphisms $\mathfrak{m} \to \mathfrak{m}(x)$. Then $F_{\gr}(\cE) := \bigcup_{x \in \cE} F_{\gr}(x)$ defines a principal fiber bundle $F_{\gr}(\cE) \to \cE$ with structure group $\aut_{\gr}(\fm)$, cf.\ \cite[Proposition~2.1]{Andreas2017}. The splitting of $D$ implies a splitting of $\fm_{-1}$, and restricting to the subgroup $G_0 \leq \aut_{\gr}(\fm)$ that preserves the splitting yields a principal subbundle $\cG_0 \to \cE$ with reduced structure group $G_0 \cong \R^{\times} \times \GL_m$, i.e., this is the {\em filtered $G_0$-structure} associated to an ODE \eqref{ODE}.

\subsection{Structure underlying the trivial ODE}\label{S:FM}

Let $n,m \geq 2$. The trivial ODE $\bu_{n+1} =\mathbf{0}$ has point symmetry Lie algebra $\g$ (see, for example, \cite[Section~2.2]{CDT2020} for explicit symmetry vector fields) with abstract structure given by
 \begin{align*}
 \g \cong \fq \ltimes V, \qquad
 \fq :=\fsl_2 \times \gl(W), \qquad
 V:= \V_n \tensor W, \qquad
 W :=\R^m,
 \end{align*}
where $\V_n$ is the unique (up to isomorphism) $\fsl_2$-irrep of dimension $n+1$, and $W$ is the standard rep of $\gl(W)$. Here, $V$ is taken to be an abelian subalgebra.

We now fix a basis for $\g$. Let $\{ w_a \}_{a=1}^m$ be the standard basis for $W$, and let $e_a{}^b$ be the $m \times m$ matrix such that $e_a{}^b w_c = \delta_c{}^b w_a$, so that $\bigl\{e_a{}^b\bigr\}_{a,b=1}^m$ spans $\gl(W)$.
Letting $\{x, y\}$ be the standard basis for $\R^2$, we identify $\V_n \cong S^n \R^2$. We obtain bases $\{ E_i \}_{i=0}^n$ on $\V_n$ and $\{ E_{i,a}\colon 0 \leq i \leq n, \, 1 \leq a \leq m \}$ on $V$ via
\begin{align*}
E_i := \dfrac{x^i y^{n-i}}{(n-i)!}, \qquad
E_{i,a} := E_i \otimes w_a.
\end{align*}
$(${\em For convenience, we define $E_i = 0$ for $i < 0$ or $i > n$. We also caution that our $E_i$ corresponds to $E_{n-i}$ in {\rm \cite[{\em Section}~2.1.2]{KT2021}.}}$)$ We complete our bases of $V$ and $\gl(W)$ to a basis of $\fg$ by introducing the standard $\fsl_2$-triple
 \begin{align*}
 \sfX := x\partial_y, \qquad
 \sfH := x\partial_x-y\partial_y, \qquad
 \sfY := y\partial_x.
 \end{align*}
Note that $\fsl_2$ commutes with $\gl(W)$, and the $\fsl_2$-actions on $\V_n$ and $V$ are naturally induced, e.g.,
 \begin{align*}
 [\sfX, E_i] = E_{i+1},\qquad
 [\sfH, E_i] = (2i-n) E_{i}, \qquad
 [\sfY, E_i] = i(n+1-i) E_{i-1}.
 \end{align*}
 In particular, $E_i$ and $E_{i,a}$ are weight vectors for the $\fsl_2$-action, i.e., eigenvectors with respect to~$\sfH$.

 Now endow $\g$ with a bi-grading as in \cite[Section~3.1]{KT2021}. Letting $\id_m := \sum_{a=1}^m e_a{}^a$, define $\sfZ_1, \sfZ_2 \in \g$ by
\begin{align}\label{Bi-grading}
\sfZ_1 := -\frac{1}{2}(\sfH + n\id_m), \qquad \sfZ_2 := -\id_m.
\end{align}
Then $\g$ decomposes into the joint eigenspaces of $\ad_{\sfZ_1}$ and $\ad_{\sfZ_2}$. We write
\begin{align*}
\g_{s,t}:= \{x \in \g\colon \sfZ_1\cdot x = sx, \, \sfZ_2 \cdot x = t x \},
\end{align*}
and refer to $s$ and $t$ as the {\em $\sfZ_1$-degree} and {\em $\sfZ_2$-degree} of $x$, respectively. The ordered pair $(s, t) \in \Z \times \Z$ is the {\em bi-grade} of $x$. It is helpful to picture $\g$ as in Figure \ref{F:Bi-grading}.

\begin{figure}[h]
\centering
\begin{tikzpicture}[scale=2,baseline=-3pt]
\draw (0.05,-1) ++(0,0) node {$\cdots$};
\draw[thick,red] (0,0) -- (1,0);
\draw[thick] (0,0) -- (-1,0);
\draw[thick] (0,0) -- (-1.5,-1);
\draw[thick] (0,0) -- (-0.5,-1);
\draw[thick] (0,0) -- (0.5,-1);
\draw[thick] (0,0) -- (1.5,-1);
\filldraw[red] (1,0) circle (0.05);
\filldraw[red] (0,0) circle (0.05);
\filldraw[black] (1.5,-1) circle (0.05);
\filldraw[black] (0.5,-1) circle (0.05);
\filldraw[black] (-0.5,-1) circle (0.05);
\filldraw[black] (-1.5,-1) circle (0.05);
\filldraw[black] (-1,0) circle (0.05);
\draw (-1,0) ++(0,-0.2) node {\tiny $(-1,0)$};
\draw (0,0) ++(0,-0.2) node {\tiny $(0,0)$};
\draw (1,0) ++(0,-0.2) node {\tiny $(1,0)$};
\draw (1.5,-1) ++(0,-0.2) node {\tiny $(0,-1)$};
\draw (0.5,-1) ++(0,-0.2) node {\tiny $(-1,-1)$};
\draw (-0.5,-1) ++(0,-0.2) node {\tiny $(-n+1,-1)$};
\draw (-1.5,-1) ++(0,-0.2) node {\tiny $(-n,-1)$};
\draw (-1,0) ++(0,0.2) node {$\sfX$};
\draw (0,0) ++(0,0.2) node {$\sfH, e_a{}^b$};
\draw (1,0) ++(0,0.2) node {$\sfY$};
\draw (-1.5,-1) ++(0,0.2) node {$E_{n,a}$};
\draw (-0.5,-1) ++(-0.1,0.2) node {$E_{n-1,a}$};
\draw (0.5,-1) ++(0.15,0.2) node {$E_{1,a}$};
\draw (1.5,-1) ++(0,0.2) node {$E_{0,a}$};
\draw [draw=black,dashed] (-0.1,-0.1) rectangle (1.1,0.1);
\end{tikzpicture}
\caption{Bi-grading on $\g$.}
\label{F:Bi-grading}
\end{figure}

Defining the {\em grading element} $\sfZ \in \fz(\g_{0,0})$, we similarly induce the structure of a $\Z$-grading on~$\g$ via
\begin{align*}
\sfZ := \sfZ_1 + \sfZ_2 = -\frac{1}{2} (\sfH + (n+2)\id_m ).
\end{align*}
(In a given representation, $\sfZ$-eigenvalues will also be referred to as {\em degrees}.) Then we have the decomposition
$\g = \g_{-n-1} \oplus \dots \oplus \g_1$, where
\begin{align*}
 &\g_1 := \g_{1,0} =\R \sfY, \\
 &\g_0 := \g_{0,0}=\R \sfH \oplus \gl_m,\\
 &\g_{-1} := \g_{-1,0} \oplus \g_{0,-1} = \R \sfX \oplus (\R E_0 \tensor W),\\
 &\g_{-i-1} := \g_{-i,-1} = \R E_i \tensor W, \qquad i = 1,\dots, n.
\end{align*}
We note that $\g_- := \g_{-n-1} \oplus \dots \oplus \g_{-1} \subset \g$ is generated by $\g_{-1}$.

We also endow $\g$ with the {\em canonical filtration} $\g^{i} := \sum_{j \ge i}\g_{j}$, which turns $\g$ into a filtered Lie algebra. Its associated graded $\gr(\g) := \bigoplus_{k \in \Z} \gr_k(\g)$, where $\gr_k(\g) = \g^k/\g^{k +1}$, is isomorphic to~$\g$ as graded Lie algebras. Using the isomorphism, we let $\gr_k\colon \g^k \to \g_k$ denote the {\em leading part}. Explicitly, if $x \in \g^k$ with $x=x_k + x_{k+1} + \cdots$, where $x_j \in \g_j$, then $\gr_k(x) := x_k$. The following notations will be convenient:
 \begin{align*}
 \fp:= \g^{0} = \fg_0 \oplus \fg_1, \qquad
 \fp_{+}:= \g^{1} = \fg_1.
 \end{align*}

At the group level, let
 \begin{align*}
 G := (\SL_2 \times \GL_m) \ltimes V, \qquad
 P := \ST_2 \times \GL_m, \qquad
 G_0 := \{ g \in P\colon \Ad_g (\g_0) \subset \g_0 \},
 \end{align*}
where $\ST_2 \subset \SL_2$ is the subgroup of lower triangular matrices. (Note that $G_0$ is isomorphic to that given in Section~\ref{S:GS}.) We also let $P_+ \subset P$ denote the connected Lie subgroup corresponding to $\fp_+ \subset \fp$. We remark that the canonical filtration on $\fg$ is $P$-invariant.

\subsection{Cartan geometries associated to ODE} \label{S:CG}

Fix $G$, $P$ and $G_0$ as above. Recall also from Section~\ref{S:GS} that all vector ODEs \eqref{ODE} can be formulated as filtered $G_0$-structures. Importantly, there is an equivalence of categories between filtered $G_0$-structures on $\cE$ (which is a wider category than that arising from ODE -- see below) and regular, normal Cartan geometries $(\cG \to \cE, \omega)$ of type $(G,P)$ \cite{CDT2020, DKM1999}. A Cartan geometry consists of a (right) principal $P$-bundle $\cG \to \cE$ endowed with a {\em Cartan connection} $\omega$, i.e., $\omega \in \Omega^1(\cG,\g)$ is a $\g$-valued 1-form on $\cG$ such that
\begin{itemize}\itemsep=0pt
\item[(a)] For any $u \in \cG$, $\omega_u\colon T_u \cG \to \fg$ is a linear isomorphism;
\item[(b)] $R_g^{\ast}\omega = \Ad_{g^{-1}} \circ\, \omega$ for any $g \in P$, i.e., $\omega$ is $P$-equivariant;
\item[(c)] $\omega(\zeta_A) = A$, where $A \in \mathfrak{p}$, where $\zeta_A$ is the fundamental vertical vector field defined by~$\zeta_A(u) := \frac{{\rm d}}{{\rm d}t}{\big|}_{t=0} u\cdot\exp(tA)$.
\end{itemize}

The {\em curvature} $K \in \Omega^2(\cG,\g)$ of the geometry is given by $K(\xi, \eta) = {\rm d}\omega(\xi, \eta) + [\omega(\xi),\omega(\eta)]$, which is $P$-equivariant and horizontal, i.e., $K(\zeta_A,\cdot) = 0$, $\forall A \in \fp$. Consequently, $K$ is determined by the $P$-equivariant {\em curvature function} $\kappa\colon \cG \to \bigwedge^2(\sfrac{\g}{\fp})^{\ast}\tensor \g$, defined by
$\kappa (u)(A, B) = K\bigl(\omega^{-1}(A),\omega^{-1}(B)\bigr)(u)$, $\forall A,B \in \g$. Letting $\omega_G$ be the Maurer--Cartan form on $G$, the {\em Klein geometry} $(G \to \sfrac{G}{P}, \omega_G)$ satisfies $K \equiv 0$ (Maurer--Cartan equation), and is the {\em flat model} for all Cartan geometries of type $(G,P)$.

 In terms of the canonical filtration $\bigl\{ \g^i \bigr\}$ on $\g$ from Section~\ref{S:FM}, $\omega$ is said to be {\em regular} if $\kappa \bigl(\g^i, \g^j\bigr) \subset \g^{i+j+1}$, $\forall i,j$. Importantly, it is known that for all filtered $G_0$-structures {\em arising from ODE}, the corresponding Cartan geometry has $\kappa$ satisfying the {\em strong regularity} condition \cite[Remark~2.3]{CDT2020}
 \begin{align}\label{E:SR}
 \kappa \bigl(\g^i, \g^j\bigr) \subset \g^{i+j+1} \cap \g^{\min(i,j)-1},\qquad
\forall i,j.
 \end{align}

To define normality, we first fix an inner product $\langle \cdot, \cdot \rangle$ on $\g$ in terms of the basis introduced in Section~\ref{S:FM}:

\begin{Definition}\label{D:metric}
Let $\langle \cdot, \cdot \rangle$ be an inner product on $\g$ such that $\bigl\{\sfX, \sfH, \sfY,e_a{}^b, E_{i,a}\bigr\}$ is an orthogonal basis for $\g$ with squared lengths of basis elements given below
\begin{align*}
\langle \sfX , \sfX \rangle = \langle \sfY, \sfY \rangle = 1, \qquad \langle \sfH, \sfH \rangle = 2, \qquad
\big\langle e_a{}^b, e_a{}^b \big\rangle = 1, \qquad \langle E_{i,a}, E_{i,a} \rangle = \frac{i!}{(n-i)!}.
\end{align*}
Then $\forall A, B \in \fq= \fsl_2 \times \gl_m$ and $\forall u, v \in V$, we have $\langle A, B \rangle = \tr\bigl(A^{\rm T}B\bigr)$ and
$\langle Au , v \rangle = \big\langle u, A^{\rm T}v \big\rangle$.
\end{Definition}

Consider $C^k(\g,\g) :=\bigwedge^k\g^{\ast} \tensor \g$ equipped with the induced canonical filtration from $\g$ and let $\partial_\g$ be the standard differential of the complex for computing Lie algebra cohomology groups~$H^k\!(\g,\!\g)$. Then, define the {\em codifferential} $\partial^\ast\colon C^k(\g, \g) \to C^{k-1}(\g, \g)$ to be the adjoint of $\partial_\g$ with respect to the induced inner product from $\g$, i.e., for each $k$ we have
$\langle \partial_\g \phi , \psi \rangle = \langle \phi, \partial^\ast \psi \rangle$ for all $\phi \in C^{k-1}(\g,\g) $ and $\psi \in C^k(\g,\g)$. By \cite[Lemma~3.2]{CDT2020}, the codifferential descends to a $P$-equivariant map $\partial^\ast\colon \bigwedge^k(\sfrac{\g}{\fp})^\ast \tensor\g \to \bigwedge^{k-1}(\sfrac{\g}{\fp})^\ast \tensor\g$. A Cartan connection $\omega$ has curvature function $\kappa$ valued in $\bigwedge^2 (\sfrac{\g}{\fp})^* \otimes \g$, and $\omega$ is said to be {\em normal} if $\partial^{\ast}\kappa = 0$. {\em In this article, we will always work with Cartan geometries of type $(G,P)$ that are normal and strongly regular.}

Since $(\partial^\ast)^2 = 0$, then the (normal) curvature $\kappa$ quotients to a $P$-equivariant func\-tion $\kappa_H$: $\cG \to \frac{\ker\partial^{\ast}}{\img \partial^{\ast}}$ called the {\em harmonic curvature}. By regularity, $\kappa_H$ is valued in the filtrand of positive degree of the $P$-module $\frac{\ker\partial^{\ast}}{\img \partial^{\ast}}$, which by \cite[Corollary~3.8]{CDT2020} is completely reducible, i.e., $P_+$ acts on it trivially, and therefore only the $G_0$-action is relevant. It is well known (see Theorem \ref{T:lh} and references therein) that $\kappa_H$ completely obstructs local flatness, i.e., $\kappa_H \equiv 0$ is equivalent to~$\kappa \equiv 0$.

Identify $\bigwedge^k (\sfrac{\g}{\fp})^* \otimes \g \cong \bigwedge^k \g_-^* \otimes \g$ as $G_0$-modules, and recall from Section~\ref{S:FM} that $\g \cong \fq \ltimes V$. Given $\phi \in C^k(\g_-,\g):=\bigwedge^k\g_-^* \tensor \g$, then we have $\phi = \sfX^* \wedge \phi_1 + \phi_2$, for $\phi_1 \in C^{k-1}(V,\g)$ and $\phi_2 \in C^k(V,\g)$, and where $\sfX^*$ is dual to $\sfX$. Denoting $\phi := \left(\begin{smallmatrix}
\phi_1\\
\phi_2
\end{smallmatrix}\right)$, then $\partial \phi $ is given by \cite[Lemma~3.4]{CDT2020}
\begin{align}\label{E:cohdef}
\partial \begin{pmatrix}
\phi_1\\
\phi_2
\end{pmatrix} = \begin{pmatrix}
-\partial_V \phi_1 + \sfX \cdot \phi_2\\
\partial_V \phi_2
\end{pmatrix},
\end{align}
where
\begin{align*}
\partial_V \phi_2(x_0,\dots,x_k) = \sum_{i=0}^k (-1)^i x_i \cdot \phi_2\bigl(x_0,\dots,\widehat{x}_i,\dots,x_k\bigr)
\end{align*}
 for $x_0,\dots,x_k \in V$, and letting $\widehat{x}_i$ denote omission of $x_i$. A direct consequence of \eqref{E:cohdef} is:
 \begin{Lemma}\label{L:X-ann}
Let $\phi\in \bigwedge^k V^*\tensor \g$. Then $\partial \phi = 0$ if and only if $\sfX \cdot \phi = 0$ and $\partial_V \phi = 0$. Moreover, if in fact $\phi \in \bigwedge^k V^* \tensor V$, then $\partial \phi = 0$ if and only if $\sfX \cdot \phi = 0$.
\end{Lemma}

 Defining $\square:=\partial \circ \partial^{\ast} + \partial^{\ast}\circ\partial \colon\bigwedge^k \g_-^* \otimes \g \to \bigwedge^k \g_-^* \otimes \g$, we then have the following $G_0$-isomorphisms,
\begin{align}\label{E:G0iso}
\bigwedge{}^{\!\!k}{\g_{-}^{\ast}\tensor \g} \cong
\rlap{$\overbrace{\phantom{\,\img \partial^{\ast} \oplus \ker\square}}^{\ker \partial^{\ast}}$}\img \partial^{\ast} \oplus
\underbrace{\ker \square \oplus \img \partial}_{\ker \partial},\qquad
\ker \square \cong \frac{\ker \partial^{\ast}}{\img \partial^{\ast}}\cong \frac{\ker \partial}{\img \partial}=:H^k(\g_{-},\g).
\end{align}

Consequently, for a regular, normal Cartan geometry, the codomain of $\kappa_H$ can be identified with the subspace $H^2_+(\g_-,\g) \subset H^2(\fg_-,\fg)$ on which the grading element $\sfZ = \sfZ_1 + \sfZ_2$ acts with positive eigenvalues. However, it should be emphasized that only part of $H^2_+(\g_-,\g)$ is in fact realizable for geometries associated to ODE \cite{DM2014, Medvedev2010}. Correspondingly, we define:

\begin{Definition}\label{D:E} The {\em effective part} $\bbE \subsetneq H^2_+(\g_{-},\g)$ is the minimal $G_{0}$-module in which $\kappa_H$ is valued, for any (strongly) regular, normal Cartan geometry of type $(G, P)$ associated to an ODE~\eqref{ODE} (for fixed $(n,m)$).
\end{Definition}

\subsection{Vector ODEs of C-class} \label{S: C-class}

We will focus on ODEs \eqref{ODE} of C-class, which have been characterized in \cite{CDT2020} using curvatures $\kappa$ of corresponding canonical Cartan connections $\omega$ described above. We define \cite[Definition~2.4]{CDT2020}:
\begin{Definition} \label{D:C-class}
An ODE \eqref{ODE} is said to be {\em of C-class} if the curvature $\kappa$ of the corresponding strongly regular, normal Cartan geometry satisfies $\kappa(\sfX,\cdot) = 0$, where $\sfX \in \g_{-1}$ was defined in Section~\ref{S:FM}.
\end{Definition}

\begin{Remark}\label{R:kC-class}
Recall from Section~\ref{S:FM} that $\g \cong \fq \ltimes V$. We remark that for a Cartan geometry corresponding to an ODE of C-class, we can identify $\kappa \in \bigwedge^2 (\g / \fq)^*\otimes \fg \cong \bigwedge^2 V^\ast \tensor \g$.
\end{Remark}

As shown in \cite{CDT2020}, the notion of C-class can be concretely reformulated in terms of fundamental invariants for vector ODEs \eqref{ODE} of order $\ge 3$ described below, which comprise the harmonic curvature of the geometry. We then have the following characterization of the C-class given in~\cite[Theorems~4.1 and~4.2]{CDT2020}:
\begin{Theorem} A vector ODE \eqref{ODE} of order $\geq 3$ is of C-class if and only if all of its generalized Wilczynski invariants vanish.
\end{Theorem}

For concreteness, we now explicitly describe the fundamental invariants for vector ODEs~\eqref{ODE} of order $n+1 \ge 3$, which consist of {\em generalized Wilczynski invariants} $\cW_r$ \cite{Doubrov2008} and {\em C-class invariants}~\cite{DM2014, Medvedev2010,Medvedev2011}:
\begin{itemize}\itemsep=0pt
\item Consider a linear vector ODE of order $n+1$:
\begin{align} \label{lODE}
\bu_{n+1} + P_n(t)\bu_n+ \dots + P_1(t)\bu_1 + P_0(t)\bu = 0,
\end{align}
where $P_j(t)$ is an $\homo(\R^m)$-valued function. Using the invertible transformations $(t,\bu) \mapsto (f(t), h(t)\bu)$ where $f\colon \R \to \R^\times$ and $h\colon \R \to \operatorname{GL}(m)$, which preserve the form of equation~\eqref{lODE}, we may normalize to $P_n = 0$ and $\tr(P_{n-1}) = 0$, i.e., {\em Laguerre--Forsyth canonical form}. Then
\begin{align*}
\Theta_{r} =
\sum_{k=1}^{r-1}(-1)^{k+1}\dfrac{(2r-k-1)!(n-r+k)!}{(r-k)!(k-1)!}P_{n-r+k}^{(k-1)}, \qquad r=2,\dots,n+1,
\end{align*}
 are fundamental invariants found by Se-ashi \cite{Se-ashi1988}, and $r$ is the degree of the invariant. For~\eqref{ODE}, the generalized Wilczynski invariants $\cW_r$ (for $r=2,\dots,n+1$) are defined as~$\Theta_r$ above evaluated at its linearization along a solution $\bu$. Formally, $\cW_r$ are obtained from~\eqref{lODE} by replacing $P_r(t)$ by the matrices $-\left(\frac{\partial f^a}{\partial u_r^b}\right)$ and the usual derivative by the total derivative $\frac{\rm d}{{\rm d}t}$ given in \eqref{E:EV}. Moreover, $\cW_r$ do not depend on the choice of solution~$\bu$, and are therefore contact invariants.
\item {\em C-class invariants} are the following:
\begin{align*}
&n \geq 2\colon\ {\bigl(\cA_2^{\tf}\bigr)}^a_{bc} =\tf\left( \frac{\partial^2{f}^a}{\partial{u_n^b}\, \partial{u_n^c}}\right), \\
&n \geq 3\colon\ {\bigl(\cA_2^{\tr}\bigr)}^a_{bc} =\tr\left( \frac{\partial^2{f}^a}{\partial{u_n^b}\, \partial{u_n^c}}\right), \\
&n = 2\colon\ (\cB_4)_{bc}= -\frac{\partial H_c^{-1}}{\partial u_{1}^b}+ \frac{\partial}{\partial u_2^b} \frac{\partial}{\partial u_2^c}{H^t}-\frac{\partial}{\partial u_2^c}\frac{d}{dt}{H^{-1}_b} \\
&\hphantom{n = 2\colon\ (\cB_4)_{bc}= }{} -\frac{\partial}{\partial u_2^c}\left(\sum_{a=1}^m H_a^{-1}\frac{\partial f^a}{\partial u_2^b}\right) + 2H^{-1}_b H^{-1}_c,
\end{align*}
where
\begin{align*}
&H_b^{-1} = \frac{1}{6(m+1)} \sum_{a=1}^m \frac{\partial^2{f}^a}{\partial{u_2^a}\, \partial{u_2^b}}, \\
&H^t = -\frac{1}{4m} \sum_{a=1}^m \left( \frac{\partial f^a}{\partial u_{1}^a}-\frac{{\rm d}}{{\rm d}t}\frac{\partial f^a}{\partial u_{2}^a}+ \frac{1}{3} \sum_{c=1}^m \frac{\partial f^a}{\partial u_2^c} \frac{\partial f^c}{\partial u_2^a}\right).
\end{align*}
\end{itemize}

 \subsection{C-class modules} \label{S:CM}

 The above fundamental invariants correspond to $G_0$-irreducible submodules in the effective part~$\bbE \subsetneq H^2_+(\g_-,\g)$ (see Definition \ref{D:E}), which we now describe.
Recall that $\g_0 \cong \sn \{\sfZ_1, \sfZ_2\} \oplus \fsl(W)$, and we have the induced action of $\sfZ_1$ and $\sfZ_2$ from \eqref{Bi-grading} on $H^2_+(\g_-,\g)$, and therefore on~$\bbE$. Note that $\sfZ_2$ acts with degrees $0$, $1$, or $2$. We define:

\begin{Definition}\label{D:C-classM}
A $G_0$-submodule $\bbU \subset \bbE \subsetneq H^2_+(\g_-,\g)$ on which $\sfZ_2$ acts with positive degree(s) is called a {\em C-class module}, and we let $\bbE_C \subsetneq \bbE$ denote the direct sum of all irreducible C-class modules. On the other hand, if $\sfZ_2$ acts on $\bbU$ with zero degree, we refer to $\bbU$ as a {\em Wilczynski module}.\looseness=-1
\end{Definition}

Any $\g_{0}$-irrep $\bbU \subset \bbE$ is determined by its bi-grade and its lowest weight $\lambda$ with respect to~$\fsl(W) \cong \fsl_m$. Such $\lambda$ can be expressed in terms of the fundamental weights $\lambda_1, \dots, \lambda_{m-1}$ of~$\fsl_m$ with respect to the Cartan subalgebra $\fh$ consisting of diagonal matrices in $\fsl_m$, and the standard choice of $m-1$ simple roots. Letting $h = \diag (h_1, \dots,h_m) \in \fh$ and $\epsilon_a\colon \fh \to \R$ the linear functional $\epsilon_a(h) = h_a$, we then have $\epsilon_1 +\dots+ \epsilon_m = 0$ and $\lambda_i = \epsilon_1 + \dots + \epsilon_i$ for $1 \leq i \leq m-1$.

Table \ref{Tab:EP} contains a summary of results for $\bbE_C \subsetneq \bbE$ for ODEs \eqref{ODE}, due to Medvedev \cite{Medvedev2010,Medvedev2011} for order $3$, and Doubrov--Medvedev \cite{DM2014} for order $\ge 4$. Using the $G_0$-isomorphisms \eqref{E:G0iso}, we identify each irreducible C-class module $\bbU \subset \bbE_C$ from Table \ref{Tab:EP} with the corresponding module in $\ker\square \subset C^2(\g_-,\g)$ consisting of harmonic 2-cochains satisfying the strong regularity con\-di\-tion~\eqref{E:SR}. (A further condition is formulated in Section~\ref{S:Ef} below.) Adopting the same notation from \cite[Table 6]{KT2021}, we let $\bbA_2$ and $\bbB_4$ denote the C-class submodules with bi-grades~$(1,1)$ and~$(2,2)$ respectively. From the respective $\sfZ_2$-degrees, and since $\kappa \in \bigwedge^2 V^* \otimes \fg$ for C-class ODE, then we deduce that we may identify
 \begin{align} \label{E:A2B4}
 \bbA_2 \subset \bigwedge{\!}^2\, V^* \otimes V, \qquad \bbB_4 \subset \bigwedge{\!}^2\, V^* \otimes \fq.
 \end{align}
Since $\bbA_2$ is not irreducible, we decompose it into (irreducible) trace and trace-free parts: $\bbA_2 = \bbA_2^{\tr} \oplus \bbA_2^{\tf}$. (The C-class invariants $\cB_4$, $\cA_2^{\tr}$, $\cA_2^{\tf}$ from Section~\ref{S: C-class} are valued in the corresponding irreducible C-class modules $\bbB_4$, $\bbA_2^{\tr}$, $\bbA_2^{\tf}$ respectively.)

 \begin{table}[h]\renewcommand{\arraystretch}{1.2}
 \centering
 $\begin{array}{|c|c|c|c|c|} \hline
 n & \text{Irred.\ C-class module}\, \bbU & \text{Bi-grade} & \fsl(W)\text{-module structure} & \fsl(W) \text{-lowest weight } \lambda\\ \hline\hline
 2 & \bbB_4 & (2,2) & S^2 W^{\ast} & -2\epsilon_1\\	
 \geq 3 & \bbA_2^{\tr} & (1,1) & W^{\ast} & -\epsilon_1\\	
 \ge 2 & \bbA_2^{\tf} & (1,1) & \bigl(S^2 W^{\ast}\tensor W\bigr)_0 & \epsilon_m - 2\epsilon_1\\\hline
 \end{array}$
 \caption{C-class modules in $\bbE_C \subsetneq \bbE \subsetneq H^2_{+}(\g_{-},\g)$ for vector ODEs of order $n+1 \geq 3$.}
 \label{Tab:EP}
 \end{table}

Since each $\bbU$ is a $\g_0$-irrep, then up to scale $\bbU$ contains a unique {\em lowest weight vector} $\Phi_\bbU$. Since $\g_0 \cong \operatorname{span}\{ \sfZ_1,\sfZ_2 \} \oplus \fsl(W)$, then being ``lowest'' means that $\Phi_\bbU$ is annihilated by all lowering operators, i.e., strictly lower triangular matrices, in $\fsl(W) \cong \fsl_m$. From Table \ref{Tab:EP}, we can give an explicit description of the annihilators $\fann(\Phi_\bbU)$, which will be needed later. Namely, if $\tilde\fp \subset \fsl_m$ is the parabolic subalgebra preserving $\Phi_{\bbU}$ up to scale, then $\fann(\Phi_\bbU) \subset \sn\{ \sfZ_1,\sfZ_2 \} \oplus \tilde\fp$. For~$a \neq c$, if $e_a{}^c \in \tilde\fp$, then $e_a{}^c \in \fann(\Phi_{\bbU})$. It suffices to consider linear combinations of $\sfZ_1$, $\sfZ_2$, and diagonal elements $\fh \subset \tilde\fp$. If $\Phi_\bbU$ has $\fsl_m$-weight $\lambda$ and $\sfZ_2$-degree $t$, then we conclude that $\fann(\Phi_{\bbU})$ is spanned by
 \begin{align} \label{dann}
 \sfZ_1 - \sfZ_2, \qquad
 h -\frac{\lambda(h)}{t}\sfZ_2, \quad h \in \fh, \qquad
 e_a{}^c \in \tilde\fp,\quad a \neq c,
 \end{align}
 where $\sfZ_1 - \sfZ_2 \in \fann(\Phi_{\bbU})$ because of the bi-grading of $\bbU$. Applying \eqref{dann} to $(\lambda,t)$ from Table~\ref{Tab:EP}, we obtain Table~\ref{Tab:fa}. Here, $\tilde\fp_1$, $\tilde\fp_{1,m-1}$ are the parabolic subalgebras in $\fsl_m$ consisting of block {\em lower} triangular matrices with diagonal blocks of sizes $1$, $m-1$ and $1$, $m-2$, $1$ respectively.

\begin{table}[h]\renewcommand{\arraystretch}{1.2}
\centering
$\begin{array}{|c|c|c|l|} \hline
n & \bbU & \dm\fann(\Phi_\bbU) & \multicolumn{1}{c|}{\text{Generators for } \fann (\Phi_{\bbU}) \subset \g_{0}} \\ \hline\hline
 2 & \bbB_4 & \multirow{2}{*}{$m^2-m+1$}&\multirow{2}{*}{$ \begin{array}{l}
\sfZ_1-\sfZ_2,\ e_a{}^c \in \tilde\fp_1,\ a \neq c, \\
 e_b{}^b-e_{b+1}{}^{b+1} + \delta_1{}^b \sfZ_2, \ 1 \le b \le m-1\\
\end{array}$}\\
\ge 3 & \bbA_2^{\tr} & & \\ \hline
 \geq 2 &\bbA_2^{\tf} & m^2-2m+3 &
 \begin{array}{l} \sfZ_1-\sfZ_2, \
 e_a{}^c \in \tilde\fp_{1,m-1}, \ a \neq c,\\
e_b{}^b-e_{b+1}{}^{b+1} + \bigl(2\delta_1{}^b + \delta_{m-1}{}^b\bigr)\sfZ_2,\ 1 \le b \le m-1\\
\end{array} \\ \hline
\end{array}$
\caption{$\fann (\Phi_\bbU) \subset \g_0$ for irreducible C-class modules $\bbU \subset \bbE$.}
\label{Tab:fa}
\end{table}

\subsection{The Doubrov--Medvedev condition} \label{S:Ef}

 We will be able to precisely identify $\bbA_2$ with the help of an additional linear condition formulated in \cite[Section~3.1, Proposition~4]{DM2014}, and which we now summarize. Consider the $\fp$-invariant subspace $F=\sn \{E_0,\dots,E_{n-1}\}\tensor W \subset V$, and define $\delta\colon \hm(F,\R\sfX) \to \hm \bigl(\bigwedge^2 F, V/F\bigr)$ by\looseness=-1
\begin{align}\label{E:delta}
(\delta B)(x,y)=(B(x)\cdot y -B(y)\cdot x)\mod F, \qquad \forall B \in \hm(F,\R\sfX).
\end{align}
We have the inclusion $\iota_F\colon F \to V$, which induces $V/F \cong W$ (as $\fp$-modules) and natural quotient $\pi_W\colon V \to V/F$. Also induced is the inclusion $\iota_{\bigwedge^2F} \colon \bigwedge^2F \to \bigwedge^2V$, from which we define $\vartheta\colon \hm \bigl(\bigwedge^2V,V\bigr) \to \hm \bigl(\bigwedge^2F,V/F\bigr) $ by $\vartheta = \pi_W \circ \iota_{\bigwedge^2 F}^*$, i.e.,
\begin{align}\label{E:alpha}
\vartheta(A) = \left.A\right|_{\bigwedge^2F} \mod F.
\end{align}
From \cite[Section~3.1, Proposition~4]{DM2014} and Remark \ref{R:kC-class}, we deduce that for a C-class ODE of order $\geq 4$, the $\bbA_2$-component $\cA_2$ of its harmonic curvature $\kappa_H$ satisfies $\vartheta(\cA_2) \in \img(\delta)$, which we refer to as the {\em Doubrov--Medvedev condition}. Correspondingly, for $n \geq 3$ we formulate the algebraic condition
\begin{align}\label{E:ALC}
\vartheta(A) \in \img(\delta), \qquad \forall A \in \bbA_2,
\end{align}
which we refer to as the {\em DM condition}. (This condition is not present for 3rd order ODE.)
 \section{Lowest weight vectors for irreducible C-class modules} \label{S:lwv}

The $\fg_0$-module structure for irreducible C-class modules $\bbU \subset \bbE$ was stated in Table \ref{Tab:EP}. While this abstract structural information proved useful in our previous study of symmetry gaps \cite{KT2021}, more precise information is needed in our current study. Namely, viewing $\bbU$ as harmonic 2-cochains via the $G_0$-equivariant identification \eqref{E:G0iso}, we may ask for concrete realizations of lowest weight vectors $\Phi_\bbU \in \bbU$ (from which a full basis of $\bbU$ may be obtained by applying raising operators). These realizations are not found in the existing literature, and our main goal in this section is to provide them. This information will provide the starting point in subsequent sections for our classification of submaximally symmetric structures.

 Given the notation introduced in Section~\ref{S:FM}, and letting $E^{i,a}$ denote the dual basis elements to $E_{i,a}$, we have:

\begin{Theorem}\label{T:lwv}
Fix $n,m \geq 2$ and an irreducible C-class module $\bbU \subset \bbE$, viewed as a $G_0$-submodule of $\ker\square \subset C^2(\fg_-,\fg)$ via \eqref{E:G0iso}. Then the unique lowest weight vector $\Phi_\bbU \in \bbU$, up to a scaling, is given in Table {\rm \ref{Tab:lwv}}.
\end{Theorem}
\begin{table}[h]\renewcommand{\arraystretch}{1.3}
\centering
$\begin{array}{|c|c|l|}\hline
n & \begin{array}{c}
\bbU\end{array}& \multicolumn{1}{c|}{\text{Lowest weight vector}\, \Phi_\bbU \in \bbU} \\\hline\hline
2 & \bbB_4 & \begin{array}{l}
E^{2,1} \wedge E^{1,1} \tensor \sfX -\frac{1}{2}E^{2,1} \wedge E^{0,1} \tensor \sfH -\frac{1}{2}E^{1,1} \wedge E^{0,1} \tensor \sfY \\ \quad + \sum_{a=1}^m \bigl(E^{2,1} \wedge E^{0,a} - E^{1,1} \wedge E^{1,a} + E^{0,1} \wedge E^{2,a} \bigr)\tensor e_a{}^1
\end{array} \\\hline
\ge 3 & \bbA_2^{\tr}&
 \begin{array}{l}
 \alpha \sum_{i=0}^n\big[\Phi^{2,i} + \bigl(\frac{n}{2}-i\bigr)\Phi^{1,i}-\frac{1}{2}i(n+1-i)\Phi^{0,i}\big] \\
\quad + \beta\sum_{i=0}^n \big[(n+1-i)\bigl(\Phi^{i,0}-\Phi^{0,i}\bigr) + \Phi^{i,1}-\Phi^{1,i}\big],\\[0.05in]
 \begin{array}{rl}
 \text{where} & \Phi^{i,j} := \sum_{a=1}^m E^{i,1}\wedge E^{j,a}\tensor E_{i+j-1,a}\\
 \text{and} & \alpha = \frac{-6(n-1)(m+1)}{mn(n+1)+6}\beta
 \end{array}
\end{array}
\\\hline
 \ge 2 & \bbA_2^{\tf}&
 \begin{array}{l}
 \sum_{j=0}^n\big[(n+1-j)\Phi^{0,j} + \Phi^{1,j}\big],\\[0.05in]
 \begin{array}{rl}
 \text{where} & \Phi^{i,j} := E^{i,1}\wedge E^{j,1}\tensor E_{i+j-1,m}
 \end{array}
 \end{array} \\ \hline
\end{array}$
\caption{Classification of lowest weight vectors $\Phi_{\bbU}$ for irreducible C-class modules $\bbU \subset \bbE$.}
\label{Tab:lwv}
\end{table}

Let us give a brief summary of the computations to follow. For $\Phi_\bbU \in \bbU$ lying in the appropriate module given in \eqref{E:A2B4}, we use the bi-grade and $\fsl(W)$-lowest weight data for $\bbU$ from Table \ref{Tab:EP} to first write a general form for $\Phi_\bbU$. (The reader should recall the bi-grades given in Section~\ref{S:FM}, e.g., $E_{i,a}$ has bi-grade $(-i,-1)$, and so $E_{i,a} \in \fg_{-i-1} \subset \fg^{-i-1}$.) We then further constrain this form by imposing additional linear conditions coming from harmonicity, strong regularity, and the DM condition \eqref{E:ALC}. (For example, since $\Phi_\bbU \in \bigwedge^2 V^* \otimes V$ in the $\bbA_2^{\tr}$, $\bbA_2^{\tf}$ cases, then $\partial \Phi_\bbU = 0$ if and only if $\sfX \cdot \Phi_\bbU = 0$ by Lemma \ref{L:X-ann}. Imposing $\sfX$-annihilation will be a detailed calculation involving the relations $\sfX \cdot E_{i,a} = E_{i+1,a}$ and $\sfX \cdot E^{i,a} = -E^{i-1,a}$.) This calculation will be involved, but we remark that in fact not all such conditions will need to be explicitly imposed:

\begin{Remark}\label{R:NS}
If $\Phi_\bbU$ can be constrained to a 1-dimensional subspace by imposing some of the conditions above, then $\Phi_\bbU$ necessarily satisfies all the remaining linear conditions (harmonicity, strong regularity, and \eqref{E:ALC}). This follows from {\em existence} of the module $\bbU \subset \bbE$ for ODE systems, which was established in \cite{DM2014, Medvedev2010}.
\end{Remark}

Let us now carry out the indicated computations and establish Theorem \ref{T:lwv} above.

\subsection[B\_4 case]{$\boldsymbol{\bbB_4}$ case} \label{S:PhiB4}
Since $\bbU := \bbB_4 \subset \bigwedge^2 V^* \otimes \fq$ has bi-grade $(2,2)$, then $\Phi_{\bbU}$ must be a linear combination of
\begin{align*}
& E^{2,a} \wedge E^{1,b} \tensor \sfX,\qquad E^{1,a} \wedge E^{0,b} \tensor \sfY,\qquad E^{2,a} \wedge E^{0,b} \tensor \sfH, \qquad E^{1,a} \wedge E^{1,b} \tensor \sfH,\\
& E^{2,a} \wedge E^{0,b} \tensor e_c{}^d, \qquad E^{1,a} \wedge E^{1,b} \tensor e_c{}^d, \qquad 1 \le a, b,c,d \le m.
\end{align*}
Since $\bbU$ has $\fsl(W)$-lowest weight $\lambda= -2\epsilon_1$ (Table \ref{Tab:EP}), then $\Phi_{\bbU}$ lies in the subspace spanned by
{\samepage\begin{align}
& E^{2,1} \wedge E^{1,1} \tensor \sfX, \qquad E^{2,1} \wedge E^{0,1} \tensor \sfH, \qquad E^{1,1} \wedge E^{0,1} \tensor \sfY,\nonumber\\
& E^{2,1} \wedge E^{0,1} \tensor e_1{}^1, \qquad E^{2,1} \wedge E^{0,1} \tensor e_a{}^a, \qquad E^{2,1} \wedge E^{0,a} \tensor e_a{}^1, \nonumber\\
& E^{2,a} \wedge E^{0,1} \tensor e_a{}^1, \qquad E^{1,a} \wedge E^{1,1} \tensor e_a{}^1.\label{E:gb}
\end{align}}
For $\fann(\Phi_{\bbU})$ from Table \ref{Tab:fa}, requiring $\fann(\Phi_{\bbU}) \cdot \Phi_{\bbU} = 0$ further constrains $\Phi_{\bbU}$ to lie in span of
\begin{align}\label{E:gB4}
&E^{2,1} \wedge E^{1,1} \tensor \sfX, \qquad E^{2,1} \wedge E^{0,1} \tensor \sfH, \qquad E^{1,1} \wedge E^{0,1} \tensor \sfY,\nonumber \\
& \sum_{a=1}^m E^{2,1} \wedge E^{0,1} \tensor e_a{}^a, \qquad \sum_{a=1}^m E^{2,1} \wedge E^{0,a} \tensor e_a{}^1, \qquad \sum_{a=1}^m E^{2,a} \wedge E^{0,1} \tensor e_a{}^1,\nonumber\\
& \sum_{a=1}^m E^{1,a} \wedge E^{1,1} \tensor e_a{}^1.
\end{align}
Let us briefly explain this. From Table \ref{Tab:fa}, $\fsl_{m-1}$ embeds into $\fann(\Phi_{\bbU})$ via $A \mapsto \diag(0,A)$, which acts trivially on the first 4 elements of \eqref{E:gb}. The remaining tensors in \eqref{E:gb} lie in a direct sum of 4 $\fsl_{m-1}$-reps equivalent to the sum of 4 copies of $\gl_{m-1}$. (Namely, consider the span of $E^{2,1} \wedge E^{0,1} \tensor e_a{}^b$, $E^{2,1} \wedge E^{0,b} \tensor e_a{}^1$, etc.) Since $\gl_{m-1} \cong \R \oplus \fsl_{m-1}$, then the aforementioned subspace contains a 4-dimensional subspace annihilated by $\fsl_{m-1}$. This is clearly spanned by the last 4 elements of \eqref{E:gB4} except taking the sum over $2 \leq a \leq m$. Finally, forcing annihilation with respect to $e_f{}^1$ for $f \geq 2$ yields \eqref{E:gB4}.

Let $\Phi_{\bbU}$ be a general linear combination of all elements of \eqref{E:gB4}, with $\mu_i$ denoting the coefficient of the $i$-th term, i.e., $\Phi_{\bbU} = \mu_1 E^{2,1} \wedge E^{1,1} \tensor \sfX + \mu_2 E^{2,1} \wedge E^{0,1} \tensor \sfH + \dots + \mu_7 \sum_{a=1}^m E^{1,a} \wedge E^{1,1} \tensor e_a{}^1$. We conclude our computation by imposing $\partial$-closedness for $\Phi_{\bbU}$ using
Lemma \ref{L:X-ann}:
 \begin{itemize}\itemsep=0pt
 \item $\sfX$-annihilation: This yields $\mu_2 = \mu_3 = -\frac{\mu_1}{2}$,
$\mu_4 = 0$ and $\mu_7 =\mu_5= -\mu_6$.
 \item $\partial_V$-closedness: $0=\partial_V\Phi_{\bbU}(E_{1,2},E_{2,1},E_{1,1}) = (\mu_5-\mu_1)E_{2,2}$, and hence $\mu_1 = \mu_5$.
 \end{itemize}
 This uniquely pins down $\Phi_{\bbU}$ (as stated in Table \ref{Tab:lwv}), up to a nonzero scaling. From Remark \ref{R:NS}, we in particular have that $\Phi_{\bbU}$ is normal and strongly regular. (The condition \eqref{E:ALC} does not apply for 3rd order ODE systems.)

\subsection[A\_2\^{}\{tr\} case]{$\boldsymbol{\bbA_2^{\tr}}$ case}
\label{S:Phitr}
This case proceeds similarly, but is more involved than the $\bbB_4$ case. In particular, more conditions are required to pin down the lowest weight vector (up to scale).

Let $n \geq 3$. Since $\bbU := \bbA_2^{\tr} \subset \bigwedge^2 V^* \otimes V$ has bi-grade $(1,1)$ and $\fsl(W)$-lowest weight $\lambda = -\epsilon_1$, then $\Phi_\bbU$ must be a linear combination of
\begin{align*}
E^{i,1} \wedge E^{j,a} \tensor E_{i+j-1,a}, \qquad 0 \le i, j \le n, \quad 1 \le i+j \le n+1,\quad 1 \le a \le m.
\end{align*}
Moreover, $\fann(\Phi_{\bbU})$ from Table \ref{Tab:fa} annihilates $\Phi_\bbU$, so $\Phi_{\bbU}$ is in fact constrained to be a linear combination of
\begin{align} \label{E:Phi_ij}
\Phi^{i,j} := \sum_{a=1}^m E^{i,1}\wedge E^{j,a}\tensor E_{i+j-1,a}.
\end{align}
Recalling our convention in Section~\ref{S:FM} that $E_k = 0$ for $k < 0$ or $k > n$, we have:

\begin{Proposition}\label{P:GS} Fix $n \geq 3$ and $m \geq 2$. Let $\bbU = \bbA_2^{\tr}$ and define $\Phi_{\bbU} = \sum_{i,j=0}^n c_{i,j}\Phi^{i,j}$ for $\Phi^{i,j}$ as in \eqref{E:Phi_ij}, where we may assume that $c_{0,0} = 0 =c_{i,j}$ for $i+j>n+1$. Since $\Phi_\bbU$ is $\partial$-closed and satisfies the strong regularity and DM conditions, then we have
\begin{align} \label{E:GS}
\begin{cases}
c_{i+1,j} + c_{i,j+1} = c_{i,j};& \text{$(XA)$: annihilation by $\sfX$;}\\
c_{i,j}{}= 0,\quad \text{for}\quad \min(i,j)\ge 3;&\text{$(SR)$: strong regularity;}\\
c_{n-1,2}= 0,\quad \text{for}\quad n\ge 4. & \text{$(DM)$: DM conditions beyond $(SR)$.}
\end{cases}
\end{align}
\end{Proposition}

\begin{proof}
 By Lemma \ref{L:X-ann}, $\partial$-closedness of $\Phi_\bbU$ is equivalent to its $\sfX$-annihilation, so using $\sfX \cdot E_{i,a} = E_{i+1,a}$ and $\sfX \cdot E^{i,a} = -E^{i-1,a}$, Leibniz rule, and re-indexing the summation, we straightforwardly obtain{\samepage
\[
0 = \sfX \cdot \Phi_{\bbU} = \sum_{i=0}^n \sum_{j=0}^n (c_{i,j}- c_{i+1,j} - c_{i,j+1})\Phi^{i,j}.
\]
This proves the first relations.}

Next, recall that $E_{k,a} \in \fg^{-k-1}$. Strong regularity \eqref{E:SR} of $\Phi^{i,j}$ forces that we have $E_{i+j-1,a} \in \fg^{\min(-i-1,-j-1)-1}$, i.e.,
\begin{align*}
&-i-j \ge \min( -i-1,-j-1)-1 \ge -\max(i,j)-2 \iff i+j \le \max(i,j)+2,
\end{align*}
or equivalently $\min(i,j) \leq 2$. All other terms are not present in the summation.

Finally, for the last relations we force \eqref{E:ALC} for $A=\Phi_\bbU$, i.e., $\vartheta(\Phi_\bbU) \in \img(\delta)$. Recall the maps~$\delta$ and~$\vartheta$ given in \eqref{E:delta} and \eqref{E:alpha}, and $F=\sn \{E_0,\dots,E_{n-1}\}\tensor W \subset V$. Modulo $F$,
 \begin{itemize}\itemsep=0pt
 \item $\displaystyle \vartheta(\Phi_\bbU) = \sum_{i,j=0}^n c_{i,j}\vartheta(\Phi^{i,j}) \equiv \sum_{i,j=0}^n c_{i,j}\Phi^{i,j}|_{\bigwedge^2F} \equiv \sum_{i=2}^{n-1} c_{i,n+1-i}\Phi^{i, n+1-i} \\ \hphantom{\vartheta(\Phi_\bbU\!)}\overset{\text{\tiny (SR)}}{\equiv} c_{2,n-1}\Phi^{2, n-1} + c_{n-1,2}\Phi^{n-1,2}.$
 \item $ \displaystyle \delta(E^{i,a} \tensor \sfX) \equiv \sum_{b=1}^m E^{i,a} \wedge E^{n-1,b} \tensor E_{n,b}$, i.e., bi-grade $(i-1,1)$ tensors for $0 \le i \le n-1$.
 \end{itemize}
 Since $\vartheta(\Phi_\bbU)$ only consists of bi-grade $(1,1)$ tensors, it suffices to examine the $(1,1)$ subspace of~$\img(\delta)$. From above, this always contains $\Phi^{2,n-1}$ (modulo $F$), but does not contain $\Phi^{n-1,2}$ when $n \geq 4$. Hence, beyond (SR), DM condition implies $\vartheta(\Phi_{\bbU}) \in \img(\delta)$, which forces $c_{n-1,2} = 0$ for $n \ge 4$.
\end{proof}

We now solve \eqref{E:GS}:
\begin{Proposition}\label{P:sol}
Fix $n \ge 3$. Then $(c_{i,j})_{0 \leq i,j \leq n}$ from Proposition {\rm \ref{P:GS}} is of the following form:
\begin{align} \label{E:sol}
\begin{split}
&c_{i,0}=\begin{cases}
(n-i+1)\beta, & 3 \le i \le n;\\
\alpha + (n-1)\beta, & i = 2;\\
\frac{n \alpha}{2} + (n-1)\beta, & i = 1;
\end{cases} \qquad
c_{0,i}=\begin{cases}
-c_{1,0}, & i =1; \\
(i-n-1)(\beta + \frac{i\alpha}{2}), & 2 \le i \le n;
\end{cases}\\
&c_{1,i}= \left(\frac{n}{2}-i\right)\alpha + \bigl(\delta_i{}^1-1\bigr)\beta, \quad 1 \le i \le n; \qquad c_{i,1} = \beta + \delta_i{}^2 \alpha, \quad 2 \le i \le n;\\
&c_{2,i} = (1-\delta_i{}^n)\alpha, \quad 2 \le i \le n; \hspace{3cm} c_{i,2}= 0, \quad 3 \le i \le n,
\end{split}
\end{align}
where $\alpha :=c_{2,n-1}$ and $\beta := c_{n,1}$, and all other coefficients are trivial.
\end{Proposition}

\begin{proof}
Since (DM) is only present for $n \geq 4$, we split our proof into two cases:
\begin{itemize}\itemsep=0pt
\item $n=3$: The system \eqref{E:GS} becomes
\begin{align*}
 &c_{1,0} + c_{0,1} = 0, \qquad
 c_{1,1} + c_{0,2} = c_{0,1},\qquad
 c_{2,0} + c_{1,1} = c_{1,0}, \\
 &c_{1,2} + c_{0,3} = c_{0,2}, \qquad
 c_{2,1} + c_{1,2} = c_{1,1}, \qquad
 c_{3,0} + c_{2,1} = c_{2,0}, \\
& c_{1,3} = c_{0,3}, \qquad
 c_{2,2} + c_{1,3} = c_{1,2}, \qquad
 c_{3,1} + c_{2,2} = c_{2,1}, \qquad
 c_{3,1} = c_{3,0}.
\end{align*}
Solving this in terms of $\alpha = c_{2,2}$ and $\beta = c_{3,1}$ gives \eqref{E:sol}.
\item $n \ge 4$:
	
{\bf Step 1}: Start with the assumed conditions $c_{0,0} = 0 = c_{i,j}$ for $i+j > n+1$, the (SR) relations, as well as the (DM) relation $c_{n-1,2} = 0$. Using (XA), determine the entries above $c_{n-1,2} = 0$ and left of $c_{2,n-1} =: \alpha$ (until the $(2,2)$-position), as shown below:
\begin{align*}
\begin{tiny}
\left( \begin{array}{ccc|cccccccc}
0 & * & * & * & * & \cdots & * & * & * & *\\
* & * & * & * & * & \cdots & * & * & * & *\\
* & * & * & * & * & \cdots & * & * & \ralpha & 0\\ \hline
* & * & * & \rz & \rz & {\color{red} \cdots} & \rz & \rz & 0 & 0\\
\vdots & \vdots & \vdots & \rz & \rz & {\color{red} \cdots} & \rz & 0 & 0 & 0\\
\vdots & \vdots & \vdots & {\color{red} \vdots} & {\color{red} \vdots} & \reflectbox{{\color{red} $\ddots$}} & \reflectbox{$\ddots$} & \vdots & \vdots & \vdots\\
* & * & * & \rz & \rz & \reflectbox{$\ddots$} & 0 & 0 & 0 & 0\\
* & * & * & \rz & 0 & \cdots & 0 & 0 & 0 & 0\\
* & * & \rz & 0 & 0 & \cdots & 0 & 0 & 0 & 0\\
* & * & 0 & 0 & 0 & \cdots & 0 & 0 & 0 & 0
\end{array} \right) \leadsto
\left( \begin{array}{ccc|cccccccc}
0 & * & * & * & * & \cdots & * & * & * & *\\
* & * & * & * & * & \cdots & * & * & * & *\\
* & * & \ralpha & \ralpha & \ralpha & \cdots & \ralpha & \ralpha & \ralpha & 0\\ \hline
* & * & \rz & 0 & 0 & \cdots & 0 & 0 & 0 & 0\\
\vdots & \vdots & {\color{red} 0} & 0 & 0 & \cdots & 0 & 0 & 0 & 0\\
\vdots & \vdots & {\color{red} \vdots} & \vdots & \vdots & \ddots & \vdots & \vdots & \vdots & \vdots\\
* & * & \rz & 0 & 0 & \cdots & 0 & 0 & 0 & 0\\
* & * & \rz & 0 & 0 & \cdots & 0 & 0 & 0 & 0\\
* & * & 0 & 0 & 0 & \cdots & 0 & 0 & 0 & 0\\
* & * & 0 & 0 & 0 & \cdots & 0 & 0 & 0 & 0
\end{array} \right).
\end{tiny}
\end{align*}
	
{\bf Step 2}: Using (XA), we have $\beta := c_{n,1} = c_{n,0}$. Use (XA) to determine the entries above~$c_{n,1}$ and left of $\gamma := c_{1,n}$ (until the $(1,1)$-position), as shown below:
\begin{align*}
\begin{tiny}
\left( \begin{array}{ccc|cccccccc}
0 & * & * & * & * & \cdots & * & * & * & *\\
* & {\color{red} (n-1)\alpha + \beta + \gamma} & {\color{red} (n-2)\alpha + \gamma} & {\color{red} (n-3)\alpha + \gamma} & {\color{red} (n-4)\alpha + \gamma} & \cdots & {\color{red} 3\alpha + \gamma} & {\color{red} 2\alpha + \gamma} & {\color{red} \alpha + \gamma} & \rgamma\\
* & {\color{red} \alpha + \beta} & \alpha & \alpha & \alpha & \cdots & \alpha & \alpha & \alpha & 0\\ \hline
* & \rbeta & 0 & 0 & 0 & \cdots & 0 & 0 & 0 & 0\\
\vdots & {\color{red} \beta} & 0 & 0 & 0 & \cdots & 0 & 0 & 0 & 0\\
\vdots & {\color{red} \vdots} & \vdots & \vdots & \vdots & \ddots & \vdots & \vdots & \vdots & \vdots\\
* & \rbeta & 0 & 0 & 0 & \cdots & 0 & 0 & 0 & 0\\
* & \rbeta & 0 & 0 & 0 & \cdots & 0 & 0 & 0 & 0\\
* & \rbeta & 0 & 0 & 0 & \cdots & 0 & 0 & 0 & 0\\
\rbeta & \rbeta & 0 & 0 & 0 & \cdots & 0 & 0 & 0 & 0
\end{array} \right).
\end{tiny}
\end{align*}
More precisely, we have
\begin{align} \label{E:c1i}
c_{1,i} =
(n-i)\alpha + \delta_i^1\beta + \gamma, \qquad 1 \leq i \leq n.
\end{align}
	
{\bf Step 3}: Using (XA), determine all entries above $c_{n,0} = \beta$. This yields
\begin{align}\label{E:ci0}
c_{i,0} = \begin{cases}
(n-i+1)\beta, & 3 \leq i \leq n;\\
\alpha + (n-1)\beta, & i=2;\\
n(\alpha + \beta) + \gamma, & i=1.
\end{cases}
\end{align}
	
{\bf Step 4}: Impose $c_{0,1} \overset{\text{\tiny (XA)}}{=} -c_{1,0} = -n(\alpha + \beta) - \gamma$. For $2 \leq i \leq n$, we have the telescoping sum
\begin{align}
&c_{0,i} - c_{0,1} = \sum_{k=2}^i (c_{0,k} - c_{0,k-1}) \overset{\text{\tiny (XA)}}{=} -\sum_{k=2}^i c_{1,k-1} \overset{\eqref{E:c1i}}{=} -\sum_{k=2}^i [(n-k+1)\alpha + \delta_{k-1}^1\beta + \gamma],\nonumber\\
&c_{0,i} = c_{0,1} -\beta - (i-1)\gamma - \alpha[ (n-1) + \dots + (n-i+1)]\nonumber\\
&\hphantom{c_{0,i} }{}= - (n+1)\beta - i\gamma - \alpha\left[ \binom{n+1}{2} - \binom{n-i+1}{2}\right]. \label{E:c0i}
\end{align}
	
{\bf Step 5}: Impose $c_{0,n} \overset{\text{\tiny (XA)}}{=} c_{1,n} = \gamma$. Solving this yields $\gamma = -\beta - \frac{n\alpha}{2}$. Substituting this into \eqref{E:c1i}, \eqref{E:ci0} and \eqref{E:c0i} then gives the stated result.\hfill $\qed$
\end{itemize}\renewcommand{\qed}{}
\end{proof}

We conclude our computation by imposing the coclosedness condition, i.e., $\partial^\ast \Phi_\bbU = 0$.

\begin{Proposition}\label{P:ne}
 Let $n \geq 3$ and $m \geq 2$. Take $\Phi_\bbU$ from Proposition {\rm \ref{P:GS}} with coefficients \eqref{E:sol}. Then
\begin{align}\label{E:alpharel}
\alpha = \frac{-6(n-1)(m+1)}{mn(n+1)+6}\beta.
\end{align}
\end{Proposition}

\begin{proof}
From Lemma \ref{L:norm}, we have
\begin{align}
&0 = \sum_{k=0}^{n-1} \frac{(n-k)(k+1)}{n(n-1)}(c_{k,2}-mc_{2,k})
+ \sum_{k=0}^{n} \frac{2k-n}{n} (mc_{1,k}-c_{k,1})\nonumber\\
&\hphantom{0 =}{} + \sum_{k=1}^n (mc_{0,k}-c_{k,0}).\label{E:nr}
\end{align}
We now substitute \eqref{E:sol} into \eqref{E:nr} and simplify. The computations are straightforward but tedious, and for the respective summations above, this leads to
 \begin{align*}
 0 &= -\frac{(n+2)\Omega}{6n(n-1)} -\frac{(n+2)\Omega}{6n} -\frac{(n+2)\Omega}{12} = -\frac{(n+2)(n+1)\Omega }{12(n-1)},
 \end{align*}
 where $\Omega := (mn(n+1)+ 6)\alpha + 6(n-1)(m+1)\beta$. This implies $\Omega = 0$, and hence the result.
\end{proof}

Combining \eqref{E:alpharel}, \eqref{E:sol}, and Proposition \ref{P:GS} uniquely determines $\Phi_{\bbU}$ (given in Table \ref{Tab:lwv} after some simplification/reorganization), up to a nonzero scaling. Note that \eqref{E:nr} (derived in Appendix \ref{S:coclosed}) was only a small part of the coclosedness condition, but using Remark \ref{R:NS}, we deduce that indeed $\partial^*\Phi_{\bbU} = 0$.

\subsection[A\_2\^{}\{tf\} case]{$\boldsymbol{\bbA_2^{\tf}}$ case} The trace-free case proceeds analogously to the trace case. Since $\bbU := \bbA_2^{\tf} \subset \bigwedge^2 V^* \otimes V$ has bi-grade $(1,1)$ and $\fsl(W)$-lowest weight $\lambda =\epsilon_m -2\epsilon_1$, then $\Phi_\bbU$ is a linear combination of
\begin{align}\label{E:btf}
\Phi^{i,j} := E^{i,1} \wedge E^{j,1} \tensor E_{i+j-1,m}, \qquad 0 \le i,j \le n,\quad 1 \le i+j \le n+1.
\end{align}
Note that $\Phi^{i,j} = -\Phi^{j,i}$ is annihilated by all elements of $\fann(\Phi_\bbU)$ given in Table \ref{Tab:fa}.

\begin{Proposition}\label{P:GS1} Fix $n,m \geq 2$.
Let $\bbU = \bbA_2^{\tf}$ and define $\Phi_\bbU = \sum_{i,j=0}^n c_{i,j} \Phi^{i,j}$ for $\Phi^{i,j}$ as in~\eqref{E:btf}, where we may assume that $c_{i,j} = -c_{j,i}$, and $c_{i,j} =0$ for $i+j>n+1$. Since $\Phi_\bbU$ is $\partial$-closed and satisfies the strong regularity and DM conditions, then we have
\begin{align}\label{E:GStf}
\begin{cases}
c_{i+1,j} + c_{i,j+1} = c_{i,j};& \text{$(XA)$: annihilation by $\sfX$;}\\
c_{i,j}{}= 0,\quad \text{for}\quad \min(i,j)\ge 3;&\text{$(SR)$: strong regularity;}\\
c_{n-1,2}=c_{2,n-1}= 0,\quad \text{for}\quad n\ge 3. & \text{$(DM')$: DM conditions beyond $(SR)$.}
\end{cases}
\end{align}
\end{Proposition}

\begin{proof}
 The proof is very similar to Proposition \ref{P:GS}, as we now explain. Recall that $\sfX \cdot E_{i,a} = E_{i+1,a}$ and $\sfX \cdot E^{i,a} = -E^{i-1,a}$, i.e., the $\sfX$-action on these basis elements is independent of the second index. Consequently, comparing \eqref{E:btf} and \eqref{E:Phi_ij}, it is immediate that $\sfX\cdot \Phi_{\bbU}=0$ yields the same conditions (XA). Strong regularity similarly does not involve the second index, and so we obtain the same conditions (SR).

Finally, let us focus on \eqref{E:ALC}. As in the proof of Proposition \ref{P:GS},
 \begin{itemize}\itemsep=0pt
 \item $\vartheta(\Phi_\bbU) \equiv \sum_{i=2}^{n-1} c_{i,n+1-i}\Phi^{i, n+1-i} \mod F$, which consists of bi-grade $(1,1)$ tensors;
 \item the bi-grade $(1,1)$ tensors in $\img(\delta)$ are spanned by $\sum_{b=1}^m E^{2,a} \wedge E^{n-1,b} \tensor E_{n,b} \mod F$.
 \end{itemize}
Since $m \geq 2$, then $\vartheta(\Phi_\bbU) \in \img(\delta)$ forces $\vartheta(\Phi_\bbU) \equiv 0$. This is automatic for $n=2$, while for $n \geq 3$, we have $c_{i,n+1-i} = 0$ for $2 \leq i \leq n-1$. Beyond (SR), we have merely $c_{2,n-1} = c_{n-1,2} = 0$.
\end{proof}

\begin{Proposition}
Fix $n \ge 2$. Then $(c_{i,j})_{0 \leq i,j \leq n}$ from Proposition {\rm \ref{P:GS1}} is of the following form:
\begin{align}\label{E:sol1}
c_{i,0} = -c_{0,i} = \begin{cases}
(n-i+1)\beta, & 2 \le i \le n\\
(n-1)\beta, & i=1
\end{cases}\qquad \text{and}\qquad
c_{i,1}= -c_{1,i} = \beta, \quad 2 \le i \le n
\end{align}
and all other coefficients are trivial.
\end{Proposition}

\begin{proof}
We split our proof into two cases:
\begin{itemize}\itemsep=0pt
\item $n=2$: the system \eqref{E:GStf} reduces to
\begin{align*}
c_{1,0} + c_{0,1}= 0, \qquad c_{2,0} = c_{1,0}, \qquad c_{0,2} = c_{0,1}, \qquad c_{1,2}=c_{0,2}, \qquad c_{2,1}=c_{2,0},
\end{align*}
and solving the system in terms of $c_{2,1}$ proves the claim.
\item $n\ge 3$: The conditions on $c_{i,j}$ in Proposition \ref{P:GS1} can be viewed as \eqref{E:GS} with ad\-di\-tion\-al\-ly~$c_{i,j} = -c_{j,i}$ (and $c_{2,n-1} = 0$ when $n=3$). Consequently, the solution to \eqref{E:GStf} can be obtained from the solution \eqref{E:sol} to \eqref{E:GS} by merely imposing $\alpha := c_{2,n-1} = 0$.\hfill $\qed$
 \end{itemize}\renewcommand{\qed}{}
 \end{proof}

 Combining \eqref{E:sol1} and Proposition \ref{P:GS1} uniquely determines $\Phi_{\bbU}$ (given in Table \ref{Tab:lwv}), up to a~nonzero scaling. As before, using Remark \ref{R:NS}, we deduce that $\partial^* \Phi_{\bbU} = 0$. This completes our proof of Theorem \ref{T:lwv}.

\section{Homogeneous structures and Cartan-theoretic descriptions} \label{S:Hom}

 Our method for proving Theorems \ref{T:I} and \ref{T:Main1} will rely on the fact that Cartan geometries (see Section~\ref{S:CG}) associated to submaximally symmetric vector ODEs are locally homogeneous (see Lemma \ref{L:LoHom}). In this section, we summarize all relevant symmetry-based facts about such geometries and their corresponding algebraic models of ODE type. We will use $G$, $P$, $G_0$ and $\g$ from Section~\ref{S:FM}, and the filtration and grading on $\g$ defined there.

\subsection{Symmetry gaps for ODE} \label{S:TPA}

An {\em infinitesimal symmetry} of a given Cartan geometry $(\cG \to \cE, \omega)$ of type $(G, P)$ is a $P$-invariant vector field $\xi \in \fX(\cG){}^P$ on $\cG$ that preserves $\omega$ under Lie differentiation, i.e., $\cL_\xi \omega = 0$. The collection of all such vector fields forms a Lie algebra, which we denote by
\begin{align*}
\finf(\cG, \omega) := \bigl\{\xi \in \fX(\cG){}^P\colon \cL_\xi \omega = 0\bigr\} \subset \fX(\cG).
\end{align*}

The submaximal symmetry dimension is
\begin{align*}
\begin{split}
\fS := \max \bigl\{\dim \finf(\cG, \omega)\colon& (\cG \to \cE, \omega) \ \text{strongly regular, normal of type} \ (G,P) \\
& \text{associated to a vector ODE $\cE$ \eqref{ODE}, with} \ \kappa_H \not \equiv 0 \bigr \}.
\end{split}
\end{align*}
Recall that $\bbE$ decomposes into $G_0$-irreducible submodules $\bbU \subset \bbE$. Analogous to $\fS$ above, we define:
\begin{align*}
\begin{split}
\fS_\bbU := \max \bigl\{\dim \finf(\cG, \omega)\colon& (\cG \to \cE, \omega) \
 \text{strongly regular, normal of type}\ (G,P) \\
& \text{associated to a vector ODE $\cE$ \eqref{ODE},} \ \text{with} \ 0 \not \equiv \img(\kappa_H) \subset \bbU \bigr\}.
\end{split}
\end{align*}

To define suitable algebraic upper bounds, we will need the following notion from \cite{KT2017}:
\begin{Definition} \label{D:TA}
Given a subspace $\fa_0 \subset \g_0$, the graded subalgebra $\fa = \pr(\g_-,\fa_0) := \fa_- \oplus \fa_0 \oplus \fa_1 \subset \g$, where $\fa_- :=\g_- = \g_{-n-1} \oplus \dots \oplus \g_{-1}$ and $\fa_1 :=\{x \in \g_1\colon [x,\g_{-1}] \subset \fa_0\}$, is called the {\em Tanaka prolongation algebra}. For $\phi$ in some $\g_0$-module, we define $\fa^\phi := \pr(\g_-,\fann(\phi))$, where $\fann (\phi) \subset \g_0$ is the annihilator of $\phi$.
\end{Definition}
Now, we define
\begin{align}\label{E:UU}
\fU := \max \bigl\{\dim \fa^\phi \colon 0 \ne \phi \in \bbE \bigr\} \qquad \text{and} \qquad \fU_\bbU := \max \bigl\{\dim \fa^\phi\colon 0 \ne \phi \in \bbU \bigr\}.
\end{align}

By \cite[Theorem~2.11]{KT2021}, we conclude that
\begin{align} \label{E:SU}
\fS \le \fU < \dim \g \qquad \text{and}\qquad \fS_\bbU \le \fU_\bbU \quad \text{for all $G_0$- irreducible modules}\quad \bbU \subset \bbE.
\end{align}
Note that $\fU = \max_{\bbU \subset \bbE} \fU_\bbU$. In fact, by \cite[Theorem~1.2]{KT2021}, in all of the {\em vector} cases we have equality:
\begin{align*}
\fS = \fU \qquad \text{and} \qquad \fS_\bbU = \fU_\bbU.
\end{align*}
Examples of some vector ODEs realizing these can be found in \cite[Tables 8 and~10]{KT2021}.

\subsection{Local homogeneity and algebraic models of ODE type}\label{S:AlgM}
Recall that a Cartan geometry $(\cG \to \cE, \omega)$ of type $(G,P)$ is said to be locally {\em homogeneous} if there exists a (left) action by a local Lie group $F$ on $\cG$ by principal bundle morphisms preserving~$\omega$ that projects onto a transitive action down on $\cE$. We then have \cite[Lemma~A.1]{KT2021}:
\begin{Lemma}\label{L:LoHom}
Fix a $G_0$-irrep $\bbU \subset \bbE$. Then any regular, normal Cartan geometry $(\cG \to \cE, \omega)$ of type $(G,P)$ with $ 0 \not\equiv \img(\kappa_H) \subset \bbU$ and $\dim (\finf(\cG, \omega)) = \fU_\bbU$ is locally homogeneous about any point $u \in \cG$ with $\kappa_H(u) \ne 0$.	
\end{Lemma}

By \cite[Section~A.1]{KT2021}, such a homogeneous Cartan geometry can be encoded {\em Cartan-the\-o\-ret\-i\-cal\-ly} by:

\begin{Definition} \label{D:AM}
An {\em algebraic model} $(\ff;\g,\fp)$ {\em of ODE type} is a Lie algebra $(\ff,[\cdot,\cdot]_\ff)$ such that:
\begin{itemize}\itemsep=0pt
\item [(i)] $\ff \subset \g$ is a filtered subspace whose associated graded $\fs:=\gr(\ff) \subset \g$ has $\fs_- = \g_-$;
\item [(ii)] $\ff^0$ inserts trivially into $\kappa(x,y) := [x,y]-[x,y]_\ff$, i.e., $\kappa(z,\cdot) = 0$ for all $z \in \ff^0$;
\item [(iii)] $\kappa$ is normal and strongly regular : $\partial^\ast \kappa = 0$ and $\kappa(\g^i, \g^j) \subset \g^{i+j+1} \cap \g^{\min (i,j)-1}$, $\forall i,j$.
\end{itemize}
\end{Definition}

Let $\cN$ denote the set of all algebraic models $(\ff;\g,\fp )$ of ODE type for fixed $(G,P)$. Then $\cN$ admits a $P$-action and is partially ordered:
\begin{itemize}\itemsep=0pt
\item[(1)] $P$-action:\, for $p \in P$ and $ \ff \in \cN$, we have $p\cdot \ff := \Ad_{p}(\ff)$. We will regard all algebraic models $(\ff;\g,\fp)$ of ODE type in the same $P$-orbit to be equivalent.
\item[(2)] Partial order relation $\le$:\, for $\ff, \widetilde{\ff} \in \cN$ regard $\ff \le \widetilde{\ff}$ if there exists a map $\ff \hookrightarrow \widetilde{\ff}$ of Lie algebras. We will focus on maximal elements in $(\cN, \le)$.
\end{itemize}

Combining \eqref{E:SU}, Lemma \ref{L:LoHom} and Definition \ref{D:AM}, we obtain the following key existence result:

\begin{Theorem} \label{T:Existance}
Fix an irreducible $G_0$-module $\bbU$ in the effective part $\bbE$ for vector ODEs \eqref{ODE} of order $\ge 3$. Then there {\em exists} an algebraic model $(\ff;\g, \fp)$ of ODE type with $ 0 \not\equiv \img(\kappa_H) \subset \bbU$ and~$ \dim \ff = \fU_\bbU = \fS_\bbU$.
\end{Theorem}

\begin{Remark}\label{R:Strategy}
 Conversely, by \cite[Lemma~4.1.4]{KT2017}, for a given algebraic model $(\ff; \g, \fp)$ of ODE type, there exists a locally homogeneous strongly regular, normal Cartan geometry $(\cG \to \cE, \omega)$ of type $(G,P)$ with $\finf(\cG, \omega)$ containing a subalgebra isomorphic to $\ff$.

 We caution that such a geometry may not arise from an ODE \eqref{ODE}. (For instance, the Doubrov--Medvedev condition must additionally hold.) Consequently, our strategy involves:
\begin{itemize}\itemsep=0pt
\item[(i)] classifying (up to the $P$-action) the corresponding algebraic models $(\ff;\g,\fp)$ of ODE type, and then
\item[(ii)] providing vector ODEs of C-class realizing these algebraic models.
\end{itemize}
\end{Remark}

A filtered linear space $\ff \subset \g$ can be described as the graph of some linear map on $\fs$ into $\g$ as follows. Let $\fs^\perp \subset \g$ be a graded subalgebra such that $\g = \fs \oplus \fs^\perp$. Then
$\ff := \bigoplus_i \sn\left\{ x + \fD(x) : x \in \fs_i \right\}$,
for some {\em unique} linear ({\em deformation}) map $\fD : \fs \to \fs^\perp$ such that $\fD(x) \in \fs^\perp \cap \g^{i+1}$ for $x \in \fs_i$.

We will use the following results from \cite[Section~A.1]{KT2021} in carrying out the classifications.

\begin{Lemma} \label{L:f0}
Let $T \in \ff^0$ and suppose that the complementary graded subspaces $\fs, \fs^\perp \subset \g$ are $\ad_T$-invariant, then the map $\fD\colon \fs \to \fs^\perp$ is $\ad_T$-invariant, i.e., $T\cdot \fD = 0 \iff \, \ad_T \circ \fD = \fD \circ \ad_T$.
\end{Lemma}

Recall from Section~\ref{S:CG} that $\kappa_H := \kappa \mod \img \partial^*$, where $\partial^*$ is the codifferential. We then have:
\begin{Proposition} \label{P:AM}
Let $(\ff; \g, \fp)$ be an algebraic model of ODE type. Then
\begin{itemize}\itemsep=0pt
\item [$(i)$] $(\ff, [\cdot,\cdot]_\ff)$ is a filtered Lie algebra.
\item [$(ii)$] $\ff^0 \cdot \kappa = 0$, i.e., $[z,\kappa(x,y)]_\ff = \kappa([z,x]_\ff,y) + \kappa(x,[z,y]_\ff)$, $\forall x,y \in \ff$, $\forall z \in \ff^0$.
\item [$(iii)$] $\fs \subset \fa^{\kappa_H}$.
\end{itemize}
\end{Proposition}

Following \cite[Section~2.2]{The2021}, we shall refer to $\ff$ as a (constrained) {\em filtered sub-deformation} of~$\fs$.

\subsection{Characterizing maximality of the Tanaka prolongation}\label{S:BFT}

Fix a $G_0$-irrep $\bbU \subset \bbE$, and recall $\fU_\bbU$ defined in \eqref{E:UU}. For vector ODEs \eqref{ODE}, $\fU_\bbU$ were computed in \cite[Section~3.4]{KT2021} using the fact that $\fU_\bbU = \dim \fa^{\Phi_\bbU}$, where $\Phi_\bbU \in \bbU$ is an extremal (lowest or highest) weight vector. For the purpose of our goal in Section~\ref{S:Embed}, we next prove that $\fU_\bbU$ is achieved precisely in this way:

\begin{Lemma} \label{L: Max}
Let $\bbU \subset \bbE$ be a $G_0$-irrep and $\Phi_\bbU \in \bbU$ be a lowest weight vector. Then, $\fU_\bbU = \dim \fa^{\Phi_\bbU}$. Moreover, if $0 \neq \phi \in \bbU$, then $\dim \fa^\phi = \fU_\bbU$ iff $[\phi]$ is contained in the $G_0$-orbit of~${[\Phi_\bbU] \in \bbP(\bbU)}$.
\end{Lemma}

\begin{proof} The proof used in \cite[Proposition~3.1.1]{KT2017} can be applied for our purposes here. (We note that the initial hypothesis of $G$ complex semisimple Lie group and $P \leq G$ a parabolic subgroup is not necessary. We use our $G_0$ here for the $G_0$ appearing there.) Over $\C$, the same proof yields the result. Over $\R$, the essential fact used in the proof is that the split-real Lie group~$\SL_m\R \subset G_0$ acts with a {\em unique closed orbit} $\cO$ (of minimal dimension) in $\bbP(\bbU)$, where $\bbU$ is an~$\SL_m\R$ irrep. (See \cite[Corollary~1]{Winther2023}.) For $\bbU \subset \bbE$ in Table \ref{Tab:EP}, the explicit orbits are
 \[\renewcommand{\arraystretch}{1.15}
 \begin{array}{|c|c|c|} \hline
 \bbU & \fsl(W)\text{-module structure} & \cO \subset \bbP(\bbU)\\ \hline\hline
 \bbB_4 & S^2 W^* & \bigl\{ \big[\eta^2\big]\colon [\eta] \in \bbP(W^*) \bigr\}\\
 \bbA_2^{\tr} & W^* & \bbP(W^*) \\
 \bbA_2^{\tf} & \bigl(S^2 W^* \otimes W\bigr)_0 & \bigl\{ \big[\eta^2 \otimes w\big]\colon [\eta] \in \bbP(W^*),\, [w] \in \bbP(W),\, \eta(w) = 0 \bigr\}\\\hline
 \end{array}
 \]
This finishes the proof.
\end{proof}

\subsection{Prolongation-rigidity}\label{S:PR}
In terms of the Tanaka prolongation algebra $\fa^\phi$ (see Definition \ref{D:TA}), we define:
\begin{Definition}\label{D:PR}
A $G_0$-module $\bbU \subset \bbE$ is said to be {\em prolongation-rigid} (PR) if $\fa_1^\phi = 0$ for {\em all} non-zero $\phi \in \bbU$.
\end{Definition}

Let $\bbU \subset \bbE$ be an irreducible C-class module (see Section~\ref{S:CM}). To study prolongation-rigidity, it suffices by Lemma \ref{L: Max} to consider the lowest weight vector $\phi = \Phi_\bbU \in \bbU$. By \cite[Lemma~3.3]{KT2021}, we have $\fa^{\Phi_\bbU}_1 = \R \sfY$ if and only if $\bbU$ has bi-grade that is a multiple of $(n,2)$. From Table \ref{Tab:EP}, the bi-grade of $\bbU$ is a multiple of $(1,1)$, so $\bbU$ is not PR if and only if $n=2$. A summary is given in Table \ref{Tab:fb}, with $\fa^{\Phi_\bbU}$ in each case, and $\fann(\Phi_\bbU)$ stated in Table \ref{Tab:fa}.

\begin{table}[ht!]\renewcommand{\arraystretch}{1.15}
\centering
$\begin{array}{|c|c|c|l|} \hline
n & \bbU & \bbU \text{ PR?} & \multicolumn{1}{c|}{\fa^{\Phi_\bbU}}\\ \hline\hline
2 & \bbB_4 &\times & \g_{-} \oplus \fann (\Phi_\bbU) \oplus \R\sfY\\
\ge 3 & \bbA_2^{\tr}&\checkmark & \g_{-} \oplus \fann (\Phi_\bbU) \\ \hline
2 & \bbA_2^{\tf} & \times & \g_{-} \oplus \fann (\Phi_\bbU) \oplus \R\sfY \\
\ge 3&\bbA_2^{\tf} &\checkmark & \g_{-} \oplus \fann (\Phi_\bbU) \\ \hline
\end{array}$
\caption{Prolongation-rigidity for irreducible C-class modules $\bbU \subset \bbE$.}\label{Tab:fb}
\vspace{-1mm}
\end{table}

\section{Embeddings of filtered sub-deformations} \label{S:Embed}

By Section~\ref{S:AlgM} above, all submaximally symmetric vector ODEs \eqref{ODE} can be encoded using algebraic models of ODE type. Consequently, proving our main results (see Theorems \ref{T:I} and~\ref{T:Main1}) boils down to classifying these corresponding algebraic models (see Theorem~\ref{T:Existance}). More precisely, in view of Lemma \ref{L: Max}, for each irreducible C-class module $\bbU \subset \bbE_C \subsetneq \bbE$ (see Definition~\ref{D:C-classM}), our goal is to classify (up to the $P$-action) all algebraic models $(\ff;\fg,\fp)$ of ODE type with $\kappa_H = \Phi_\bbU \in \bbU$, where $\Phi_\bbU$ is the lowest weight vector from Table \ref{Tab:lwv}, and $\dim \ff = \fS_\bbU$.

In this section, we classify all possible (filtered) linear embeddings $\ff \subset \fg$ for such $(\ff;\fg,\fp)$. The possibilities for curvature $\kappa$ of $(\ff;\fg,\fp)$ are then classified in Section~\ref{S:curvature}. Recall the canonical filtration and the grading structure on $\g$ from Section~\ref{S:FM}. Having computed graded subalgebras $\fa^{\Phi_\bbU} \subset \g$ in Table \ref{Tab:fb}, we next classify, up to the $P$-action, possible filtered linear subspaces $\ff \subset \g$ for algebraic models $(\ff;\fg,\fp)$ satisfying $\gr(\ff) = \fa^{\Phi_\bbU}$:

\begin{Proposition} \label{P:DR}
Fix an irreducible C-class module $\bbU = \bbB_4$, $\bbA_2^{\tf}$ or $\bbA_2^{\tr}$ in $\bbE_C \subsetneq \bbE$, viewed as a $G_0$-submodule of $\ker\square \subset C^2(\fg_-,\fg)$ via \eqref{E:G0iso}, and consider an algebraic model $(\ff;\fg,\fp)$ of ODE type with $ \kappa_H = \Phi_{\bbU}\in \bbU$, for $\Phi_\bbU$ from Table {\rm \ref{Tab:lwv}}, and $\dim \ff = \fS_{\bbU}$. Then using the $P$-action, $\ff \mapsto \Ad_p \ff$, we may normalize to
\begin{itemize}\itemsep=0pt
\item[$(a)$] $\ff = \fa^{\Phi_{\bbU}}$ when $\bbU = \bbB_4$ or $\bbA_2^{\tf};$
\item [$(b)$] $\ff = \sn\{E_{n,a}, \dots, E_{2,a},\, E_{1,1}+(n-2)\zeta\sfZ_1, \, E_{1,b}, E_{0,1}+\zeta\sfY,\, E_{0,b}, \, \sfX\colon \zeta \in \R,\ 1\le a \le m, \allowbreak 2 \le b \le m \} \oplus \fann (\Phi_{\bbA_2^{\tr}})$ when $\bbU = \bbA_2^{\tr}$.
\end{itemize}	
\end{Proposition}

 \begin{proof}
 Since $\bbU$ has a bi-grade that is a multiple of $(1,1)$ (see Table \ref{Tab:EP}), then $ T:=\sfZ_1 -\sfZ_2 \in \fann (\Phi_{\bbU})=\fa^{\Phi_\bbU}_0$. We note that $\fa^{\Phi_\bbU} \subset \g$ is a graded subalgebra, and denote by $\widehat{T}$ the element in $\ff^0$ with the {\em leading part} $T$, i.e., $\gr_0\bigl(\widehat{T}\bigr) = T$. Since $\g_1 = \R \sfY$ and $\g_i = 0$ for all $i \ge 2$ (see Section~\ref{S:FM}), then necessarily $\widehat{T}= T + s \sfY \in \ff^0$. We claim that without loss of generality, i.e., using the $P$-action (for $P$ defined in Section~\ref{S:FM}), we may assume that $T \in \ff^0$. This is immediate when $\bbU$ is {\em not} PR, since $\sfY \in \ff^0$ and therefore $T= \widehat{T}-s\sfY$ is a linear combination of $\widehat{T}$ and $\sfY$. Otherwise, when $\bbU$ is PR $\bigl(\sfY \not\in \ff^0\bigr)$, using the $P_{+}$-action and $[\sfY, T] = -\sfY$, we have
 \begin{align*}
\Ad_{\exp(t\sfY)}\bigl(\widehat{T}\bigr) = \exp(\ad_{t\sfY})\bigl(\widehat{T}\bigr)= \widehat{T} + \big[t\sfY,\widehat{T}\big]+ \frac{1}{2!}\big[t\sfY,\big[t\sfY, \widehat{T}\big]\big] + \cdots= \sfZ_1 - \sfZ_2 + (s-t)\sfY,
 \end{align*}
 then choosing $t=s$ normalizes the right-hand side to $T$. So, relabeling the left hand side by $\widehat{T}$, gives $T = \widehat{T} \in \ff^0$.

By Definition \ref{D:AM}\,(ii), we have
$\kappa(T, z) = 0$, i.e., $[T, z]_\ff = [T, z]$, $\forall z \in \ff$.
Then, by exploiting the semi-simplicity of $\ad_T$, we next determine the remaining basis elements $\widehat{x} \in \ff^i$ with the leading parts $x \in \fa^{\Phi_\bbU}_i$, i.e., $\gr_i(\widehat{x}) = x$.

We first consider $\widehat{x} \in \ff^0$. We claim that without loss of generality, as it was for $\widehat{T} \in \ff^0$ above, we may assume that $x \in \ff^0$ for all $x \in \fa^{\Phi_\bbU}_0$. We let $\widehat{x} = x + c_x \sfY$. Then, for $\bbU$ that is {\em not} PR the claim holds, since $\sfY \in \ff^0$, and so $x = \widehat{x}-c_x\sfY \in \ff^0$. In order to give the argument for the case when $\bbU$ is PR, we recall that $[T, \sfY] = \sfY$ and $[T, x] = 0$ for all $x \in \fa_0^{\Phi_\bbU}$. So, $[T, \widehat{x}]_\ff =[T, \widehat{x}] = c_x \sfY \in \ff^0$. Now, since $(\ff, [\cdot, \cdot]_\ff)$ is a Lie algebra and $\sfY \not \in \ff^0$, then the closure condition $[T, \widehat{x}] \in \ff^0$ implies that $c_x = 0$. So, $x = \widehat{x} \in \ff^0$.

Next, we similarly consider $\widehat{x} \in \ff^i$ for $i <0$. Recall that by Definition \ref{D:TA} for these cases we have $\fa^{\Phi_\bbU}_i = \g_i$, for $\g_i$ as was defined in Section~\ref{S:FM}. In view of Lemma \ref{L:f0}, we fix $\ad_T$-invariant subspaces $\fs^\perp \subset \g$ in Table \ref{Tab:sperp} such that $\g = \fs \oplus \fs^\perp$, where $\fs := \fa^{\Phi_\bbU}$, and define the deformation map $\fD\colon \fs \to \fs^\perp$ (see Section~\ref{S:AlgM}). Let $E^{i,a}$ and $\sfX^\ast$ denote the dual basis elements to~$E_{i,a}$ and~$\sfX$, respectively, and recall bi-grades for the basis elements from Figure \ref{F:Bi-grading}. Since, for $0 \le i \le n$, $1 \le a,b,c \le m$, the eigenvalues of $\ad_T$ on
\begin{align}\label{jac}
&E^{i,a} \tensor \sfZ_1, \qquad E^{i,a} \tensor e_c{}^b,\qquad E^{i,a} \tensor \sfY, \qquad \sfX^\ast \tensor \sfZ_1,\qquad \sfX^\ast \tensor e_c{}^b,\qquad \sfX^\ast \tensor \sfY
\end{align}
are $i-1$, $i-1$, $i$, $1$, $1$, $2$ respectively, then we have zero eigenvalues only when $i=0$ or $1$. Then $T \cdot \fD = 0$ (see Lemma \ref{L:f0}) implies $\sfX =\widehat{\sfX} \in \ff$ and $E_{i,a} = \widehat{E}_{i,a} \in \ff $ for all $i$ except possibly when $i=0$ or $1$.
\begin{table}[ht!]\renewcommand{\arraystretch}{1.2}
\centering
$\begin{array}{|c|c|c|c|} \hline
n & \text{Irreducible C-class module } \bbU & \text{Generators for}\,\fs^\perp \subset \g & \text{Ranges}\\ \hline\hline
2 & \bbB_4 & \sfZ_1, \ e_1{}^b & 2 \le b \le m \\
\ge 3 & \bbA_2^{\tr}&\sfZ_1, \ e_1{}^b, \ \sfY & 2 \le b \le m \\
2 & \bbA_2^{\tf} &
\sfZ_1, \ e_1{}^b, \ e_d{}^m & 2 \le b \le m,\ 2 \le d \le m-1 \\
\ge 3 & \bbA_2^{\tf}&\sfZ_1,\ e_1{}^b, \ e_d{}^m, \ \sfY & 2 \le b \le m,\ 2 \le d \le m-1\\ \hline
\end{array}$
\caption{$\ad_T$-invariant subspace $\fs^\perp \subset \g$ complementary to $\fs=\fa^{\Phi_{\bbU}}$.}
\label{Tab:sperp}
\end{table}	

 Now, consider the above exceptional cases. Based on the eigenvalues for $\ad_T$ given in \eqref{jac}, we must have
 \begin{align*}
 \widehat{E}_{0,a} = E_{0,a} + \lambda_a\sfY.
 \end{align*}
 Since $\bbU$ is a C-class module, then $\kappa(\sfX, \cdot) = 0$ (see Definition \ref{D:C-class}), which implies $[\sfX, \cdot]_\ff = [\sfX, \cdot]$. Recall that $(\ff, [\cdot, \cdot]_\ff)$ is a Lie algebra and $[\sfX, E_{0,a}] = E_{1,a}$. Then for $\bbU$ that is
\begin{itemize}\itemsep=0pt
\item[(a)] {\em not} PR \big($\bbU = \bbB_4$ or $\bbA_2^{\tf}$ when $n=2$\big): we have that $E_{0,a} = \widehat{E}_{0,a}- \lambda_a\sfY \in \ff$, since $\sfY \in \ff$ and $\widehat{E}_{0,a} \in \ff$. Since $\sfX \in \ff$, then $[\sfX, E_{0,a}]_\ff = [\sfX, E_{0,a}] = E_{1,a} \in \ff$. Hence, $E_{1,a} = \widehat{E}_{1,a}$ and so $\ff = \fa^{\Phi_{\bbU}}$.
\item [(b)] PR \big($\bbU =\bbA_2^{\tf}$ or $\bbA_2^{\tr}$ when $n\ge 3$\big): Recall from our discussion above that for {\em any} $\widehat{x} \in \ff^0$ with the leading part $x \in \fa^{\Phi_\bbU}_0$, we may assume without loss of generality that $x \in \ff^0$. Hence, Table \ref{Tab:fa} yields
 \begin{align*} 
 q &:= e_1{}^1 - e_2{}^2 + \sfZ_2 \in \ff^0 \quad \text{for} \quad \bbA_2^{\tr},\\
 p &:= e_1{}^1 - e_2{}^2 + (2 + \delta_{m-1}{}^1)\sfZ_2 \in \ff^0\, \quad \text{for} \quad \bbA_2^{\tf}.
 \end{align*}
Recall that by Definition \ref{D:AM}\,(ii) we have $\kappa\bigl(\ff^0,\cdot\bigr) = 0$, which implies $[z, \cdot]_\ff = [z,\cdot]$ for all $z \in \ff^0$. Now, since both $p$ and $q$ commute with $\sfY$, $[\sfZ_2, E_{i,a}] = -E_{i,a}$ and $\big[e_1{}^1 - e_2{}^2, E_{i,a}\big] = \bigl(\delta_a{}^1-\delta_a{}^2\bigr) E_{i,a}$ (see Section~\ref{S:FM}), then for
\begin{itemize}\itemsep=0pt
\item[(i)] $\bbU = \bbA_2^{\tf}$: we have $\big[p,\widehat{E}_{0,a}\big]_\ff = \big[p,\widehat{E}_{0,a}\big] = \bigl(\delta_a{}^1-\delta_a{}^2 - \delta_{m-1}{}^1-2\bigr)E_{0,a}$. So, the closure condition $\big[p,\widehat{E}_{0,a}\big] \in \ff$ forces $\lambda_a = 0$, i.e., $E_{0,a} = \widehat{E}_{0,a} \in \ff$. Then $[\sfX, E_{0,a} ]= E_{1,a} \in \ff$, which implies $E_{1,a} =\widehat{E}_{1,a} \in \ff$. Hence, $\ff = \fa^{\Phi_{\bbA_2^{\tf}}}$. This completes the proof for (a).
\item [(ii)] $\bbU = \bbA_2^{\tr}$: we have $\big[q, \widehat{E}_{0,a}\big]_\ff = \big[q, \widehat{E}_{0,a}\big] = \bigl(\delta^1_a - \delta^2_a - 1\bigr) E_{0,a}$. So, $\big[q, \widehat{E}_{0,a}\big] \in \ff$ implies that $\lambda_a = 0$ {\em except} for $a=1$, so for these cases we have $E_{0,a} = \widehat{E}_{0,a} \in \ff$. Consequently, $[\sfX, E_{0,a}] = E_{1,a} \in \ff$ implies that $E_{1,a} = \widehat{E}_{1,a} \in \ff$ except for $a=1$.

Finally, we consider the case when $a =1$. Based on the eigenvalues for $\ad_T$ given in~\eqref{jac} and setting $\lambda_1 = \zeta$, we necessarily have
\begin{align*}
&\widehat{E}_{0,1} = E_{0,1} + \zeta \sfY \qquad \text{and} \qquad \widehat{E}_{1,1} = E_{1,1} + \beta \sfZ_1 +\sum_{b=2}^m\alpha_b e_1{}^b .
 \end{align*}
Since from \eqref{Bi-grading} we have $\sfH = n\sfZ_2-2\sfZ_1$ and $T = \sfZ_1- \sfZ_2$, then we have
\begin{align*}
\big[\sfX, \widehat{E}_{0,1}\big]_\ff &= [\sfX, E_{0,1} + \zeta \sfY] = E_{1,1} + \zeta \sfH = E_{1,1} +\zeta (n\sfZ_2-2\sfZ_1)\\
&= E_{1,1} + \zeta (n-2)\sfZ_1 -\zeta nT.
\end{align*}
So, the closure condition $[\sfX, \widehat{E}_{0,1}]_\ff \in \ff$ holds only if $\widehat{E}_{1,1} = E_{1,1} + \zeta (n-2)\sfZ_1 -\zeta nT$. This implies $\beta = (n-2)\zeta$ and $\alpha_b = 0$ for {\em all} $b$, which proves (b) and concludes our proof.\hfill $\qed$
\end{itemize}
\end{itemize} \renewcommand{\qed}{}
\end{proof}

We have the following result for the curvature $\kappa$ of such an algebraic model. This result is essential for our study of curvatures in Section~\ref{S:curvature}.

\begin{Corollary}\label{C:C-class}
Fix an irreducible C-class module $\bbU \subset \bbE_C \subsetneq \bbE$, viewed as a $G_0$-submodule of $\ker\square \subset C^2(\fg_-,\fg)$ via \eqref{E:G0iso}, and consider an algebraic model $(\ff;\fg,\fp)$ of ODE type with $ \kappa_H = \Phi_{\bbU}\in \bbU$, for $\Phi_\bbU$ from Table {\rm \ref{Tab:lwv}}, and $\dim \ff = \fS_{\bbU}$, normalized according to Proposition {\rm \ref{P:DR}}. Then $\sfX \cdot \kappa = 0$.
\end{Corollary}

\begin{proof}
Since $\bbU$ is a C-class module, then, $\kappa \in \bigwedge^2V^*\tensor \g$ (see Remark \ref{R:kC-class}) for $V$ from Section~\ref{S:FM}. So, for $\sfX \in \ff$ we have $[\sfX, z]_\ff = [\sfX, z]$ for all $z \in \ff$. Then, as a consequence of the Jacobi identity we get the claim as follows:
\begin{align*}
(\sfX \cdot \kappa)(x,y) &= [\sfX, \kappa(x,y)]-\kappa([\sfX, x], y)-\kappa(x, [\sfX , y])\\
&= [\sfX, [x,y]]-\underbrace{[\sfX, [x,y]_\ff]}_{[\sfX, [x,y]_\ff]_\ff} + [\underbrace{[\sfX, x]}_{[\sfX,x]_\ff},y]_\ff-[[\sfX, x],y] + [x,\underbrace{[\sfX, y]}_{[\sfX,y]_\ff}]_\ff-[x,[\sfX, y]]=0.\!\!\!\!\!\!\tag*{\qed}
\end{align*} \renewcommand{\qed}{}
\end{proof}

 \section[Classification of submaximally symmetric vector ODEs of C-class]{Classification of submaximally symmetric vector ODEs\\ of C-class}
\label{S:curvature}

In this section, we classify (up to the $P$-action) all algebraic models of ODE type for submaximally symmetric vector ODEs \eqref{ODE} of C-class (see the introduction to Section~\ref{S:Embed}) and establish Theorems \ref{T:I} and \ref{T:Main1}.

\subsection{Algebraic curvature constraints} \label{S:k}

For the algebraic models $(\ff;\g,\fp)$ whose possible filtered linear subspaces $\ff \subset \g$ have been classified in Proposition \ref{P:DR}, we classify their possible curvatures $\kappa$ below.

\begin{Proposition}\label{P:k}
Fix an irreducible C-class module $\bbU = \bbB_4$, $\bbA_2^{\tf}$ or $\bbA_2^{\tr}$ in $\bbE_C \subsetneq \bbE$, viewed as a $G_0$-submodule of $\ker\square \subset C^2(\fg_-,\fg)$ via \eqref{E:G0iso}, and consider an algebraic model $(\ff; \g, \fp)$ of ODE type with $\kappa_H = \Phi_{\bbU} \in \bbU$, for $\Phi_\bbU$ from Table {\rm \ref{Tab:lwv}}, and $\dim \ff = \fS_{\bbU}$, normalized according to Proposition {\rm \ref{P:DR}}. Then $\kappa$ is
\begin{itemize}\itemsep=0pt
\item [$(a)$] $\bbU = \bbB_4\colon \kappa = \pm\Phi_{\bbU}$ $($over $\C$, we can take $\kappa = \Phi_{\bbU})$;
\item[$(b)$] $\bbU = \bbA_2^{\tf}\colon\kappa = \Phi_{\bbU}$;
\item[$(c)$] $\bbU = \bbA_2^{\tr}\colon\kappa = \Phi_{\bbU} + \kappa_4$, where
\begin{align}
\kappa_4 ={}& \mu_1 E^{3,1} \wedge E^{0,1} \tensor \sfX + \mu_2 E^{2,1} \wedge E^{1,1} \tensor \sfX \nonumber \\
& -\frac{\mu_1 + \mu_2}{2} \left( E^{2,1} \wedge E^{0,1} \tensor \sfH + E^{1,1} \wedge E^{0,1} \tensor \sfY \right)\nonumber\\
& + \mu_3\sum_{a=1}^m \left(E^{2,1} \wedge E^{0,a} - E^{2,a} \wedge E^{0,1}+ E^{1,a} \wedge E^{1,1} \right) \tensor e_a{}^1\label{E:kappa4}
\end{align}
for some $\mu_1,\mu_2,\mu_3 \in \R$.
\end{itemize}
\end{Proposition}

\begin{proof}
The majority of the proof will consist of evaluating the annihilation conditions $\ff^0 \cdot \kappa = 0$.

Recall from Table \ref{Tab:EP} that $\bbU$ has bi-grade either $(1,1)$ or $(2,2)$, so $\sfZ_1 - \sfZ_2 \in \fann(\Phi_\bbU)$, which is contained in $\ff^0$ by Proposition \ref{P:DR}. Hence, $(\sfZ_1 - \sfZ_2) \cdot \kappa = 0$ implies that $\kappa$ is the sum of terms with bi-grades that are multiples of $(1,1)$. But $\kappa$ is regular, lies in $\bigwedge^2(\g/\fp)^\ast \tensor \g$, and~$\sfZ_2$ acts on the latter with eigenvalues $0$, $1$ or $2$. Thus, the terms in $\kappa$ can only have bi-grades~$(1,1)$ or~$(2,2)$.
By Theorem \ref{T:lh}, $\kappa_H$ can be identified with the lowest $\sfZ$-degree component of $\kappa$. Moreover, $\kappa_H$~is a~nonzero multiple of $\Phi_\bbU$. Using the $G_0$-action by $\exp(\sfZ t)$, where~$\sfZ \in \fz(\g_0)$ is the grading element (see Section~\ref{S:FM}), this multiple can be re-scaled to~$\pm 1$. For $\bbU = \bbA_2^{\tf}$ or $\bbA_2^{\tr}$, we can further normalize this multiple to $+1$. (Use the diagonal elements in~${g = \diag(a_1,\dots ,a_m) \in \GL_m \subset G_0}$, i.e., $g \cdot \Phi^{i,j} = \frac{1}{a_1} \Phi^{i,j}$ for $\Phi^{i,j}$ from \eqref{E:Phi_ij}, while
${g \cdot \Phi^{i,j} = \frac{a_m}{(a_1)^2} \Phi^{i,j}}$ for $\Phi^{i,j}$ from \eqref{E:btf}.)
Summarizing, we have
\begin{align*}
\kappa =
\begin{cases}
\pm \Phi_\bbU, & \text{when} \quad \bbU = \bbB_4,\\
\Phi_\bbU + \kappa_4, & \text{when} \quad \bbU = \bbA_2^{\tf} \quad \text{or}\quad \bbA_2^{\tr},
\end{cases}
\end{align*}
where $\kappa_4$ is the bi-grade $(2,2)$ component of $\kappa$. The $\bbB_4$ case is complete, and we turn to the remaining cases.
	
Since $\bbU$ is a C-class module, then by Remark \ref{R:kC-class} we have $\kappa \in \bigwedge^2 V^* \otimes \fg$ in the notation of Section~\ref{S:FM}. Recall $\fg = \fq \ltimes V$, and $\fq$ and $V$ have $\sfZ_2$-degrees $0$ and $-1$ respectively (see Figure~\ref{F:Bi-grading}). In particular, $\Phi_\bbU \in \bigwedge^2 V^\ast\tensor V$ and $\kappa_4 \in \bigwedge^2 V^\ast \tensor \fq$. More precisely, since $\kappa_4$ has bi-grade $(2,2)$, then in terms of the dual basis elements $E^{i,a}$ to $E_{i,a}$, having bi-grades $(i,1)$ and $(-i,-1)$ respectively, $\kappa_4$ must lie in the subspace $K_4 \subset \bigwedge^2 V^\ast \tensor \fq$ spanned by
\begin{align}
& E^{1,a} \wedge E^{0,b} \tensor \sfY, \qquad E^{3,a} \wedge E^{0,b} \tensor \sfX, \qquad E^{2,a} \wedge E^{1,b} \tensor \sfX,\nonumber\\
& E^{2,a} \wedge E^{0,b} \tensor \sfH, \qquad E^{1,a} \wedge E^{1,b} \tensor \sfH, \qquad E^{2,a} \wedge E^{0,b} \tensor e_c{}^d, \qquad E^{1,a} \wedge E^{1,b} \tensor e_c{}^d,\label{K4}
\end{align}
where $1 \le a, b,c,d \le m$. We will further constrain $\kappa_4$ as follows. Using Proposition \ref{P:DR}, we have~$\fann(\Phi_\bbU) \subset \ff^0$. Such elements annihilate both $\Phi_\bbU$ and $\kappa$, and so
\begin{align*}
z \cdot \kappa_4 = 0, \qquad \forall z \in \fann(\Phi_\bbU).
\end{align*}
	
Let us use these to find more explicit conditions on $\kappa_4$.
\begin{itemize}\itemsep=0pt
\item [(1)] $\bbU = \bbA_2^{\tf}$: from Table {\rm \ref{Tab:fa}}, we have $p:=e_1{}^1 - e_2{}^2 + \bigl(2 + \delta_{m-1}{}^1\bigr)\sfZ_2 \in \fann(\Phi_\bbU)$. From $0 = p\cdot \kappa_4$ and $\sfZ_2 \cdot \kappa_4 = 2\kappa_4$, we have that $\kappa_4$ has eigenvalue $\lambda = -2\bigl(2 + \delta_{m-1}{}^1\bigr)$ for $e_1{}^1 - e_2{}^2$. We conclude that $\kappa_4 = 0$ $($hence $\kappa = \Phi_\bbU)$ from the following considerations:
\begin{itemize}\itemsep=0pt
\item[$(i)$] $m=2$: We have $\lambda = -6$. Noting that $e_1{}^1 - e_2{}^2$ commutes with $\{\sfX, \sfH, \sfY\}$, and
\begin{align*}
&\bigl(e_1{}^1 - e_2{}^2\bigr)\cdot E^{i,a} = \bigl(\delta_a{}^2-\delta_a{}^1\bigr)E^{i,a}, \\
&\bigl(e_1{}^1 - e_2{}^2\bigr)\cdot e_a{}^b = \delta_a{}^1e_1{}^b-\delta_1{}^b e_a{}^1-\delta_a{}^2e_2{}^b + \delta_2{}^b e_a{}^2.
\end{align*}
From \eqref{K4}, we conclude that the eigenvalues of $e_1{}^1 - e_2{}^2$ in $K_4$ lie between~$-4$ and~$4$. Since $-6$ is not an eigenvalue, then $\kappa_4 = 0$.
			
\item[$(ii)$] $m\geq 3$: We have $\lambda = -4$. Proceeding as in $(a)$, we observe that $K_4$ has $-4$-eigenspace for $e_1{}^1 - e_2{}^2$ spanned by $E^{2,1} \wedge E^{0,1} \tensor e_2{}^1$. But from Table {\rm \ref{Tab:fa}}, we also have $e_{m-1}{}^{m-1} -e_m{}^m + \sfZ_2 \in \fann(\Phi_\bbU)$, which must similarly annihilate $\kappa$ and $\kappa_4$. But its eigenvalue on $E^{2,1} \wedge E^{0,1} \tensor e_2{}^1$ is $\delta_{m-1}{}^2 +2$, which is nonzero, so $\kappa_4 = 0$ follows.
\end{itemize}

\item [$(2)$] $\bbU = \bbA_2^{\tr}$: from Table \ref{Tab:fa}, we have $q_d := e_d{}^d - e_{d+1}{}^{d+1} + \delta_1{}^d \sfZ_2 \in \fann(\Phi_{\bbU})$, so $0 = q_d \cdot \kappa_4$ for $1 \leq d \leq m-1$. Letting $\fh \subset \fsl_m$ denote the standard Cartan subalgebra consisting of diagonal trace-free matrices, and $\epsilon_a \in \fh^*$ the standard weights for $\fh$, we observe
 \begin{itemize}\itemsep=0pt
 \item[(i)] $0 = q_d \cdot \kappa_4$ for $1 \leq d \leq m-1$ is equivalent to $\kappa_4$ having weight $-2\epsilon_1$,
 \item[(ii)] the first five elements of \eqref{K4} have weight $-\epsilon_a - \epsilon_b$,
 \item[(iii)] the last two elements of \eqref{K4} have weight $-\epsilon_a - \epsilon_b + \epsilon_c - \epsilon_d$.
 \end{itemize}
 Matching these weights with $-2\epsilon_1$, we deduce that $\kappa_4$ lies in the span of the following:
\begin{align*}
&E^{3,1} \wedge E^{0,1} \tensor \sfX, \qquad E^{2,1} \wedge E^{1,1} \tensor \sfX, \qquad E^{2,1} \wedge E^{0,1} \tensor \sfH,\\
&E^{1,1}\wedge E^{0,1} \tensor \sfY, \qquad E^{2,1} \wedge E^{0,1} \tensor e_1{}^1, \qquad E^{2,1} \wedge E^{0,1} \tensor e_a{}^a, \\
&E^{2,1} \wedge E^{0,a} \tensor e_a{}^1, \qquad E^{2,a} \wedge E^{0,1} \tensor e_a{}^1, \qquad E^{1,a} \wedge E^{1,1} \tensor e_a{}^1,
\end{align*}
where $2 \leq a \leq m$. Similarly as in Section~\ref{S:PhiB4}, we conclude that imposing annihilation by all of $\fann(\Phi_{\bbA_2^{\tr}}) \subset \ff^0$ forces $\kappa_4$ to lie in the subspace spanned by
\begin{align}
&E^{3,1} \wedge E^{0,1} \tensor \sfX, \qquad E^{2,1} \wedge E^{1,1} \tensor \sfX, \qquad E^{2,1} \wedge E^{0,1} \tensor \sfH, \qquad E^{1,1} \wedge E^{0,1} \tensor \sfY, \nonumber\\
& \sum_{a=1}^m E^{2,1} \wedge E^{0,1} \tensor e_a{}^a, \qquad \sum_{a=1}^m E^{2,1} \wedge E^{0,a} \tensor e_a{}^1, \qquad \sum_{a=1}^m E^{2,a} \wedge E^{0,1} \tensor e_a{}^1, \nonumber\\
& \sum_{a=1}^m E^{1,a} \wedge E^{1,1} \tensor e_a{}^1.\label{E:offdiag}
\end{align}

Finally, we complete the proof by imposing $\sfX \cdot \kappa = 0$ (see Corollary \ref{C:C-class}). Since $\sfX \cdot \Phi_{\bbA_2^{\tr}} = 0$ (see Lemma \ref{L:X-ann}), then $\sfX \cdot \kappa = 0$ implies that $\sfX \cdot \kappa_4 = 0$. Now let $\kappa_4$ be a general linear combination of all elements of \eqref{E:offdiag}, i.e., $\kappa_4 = \nu_1 E^{3,1} \wedge E^{0,1} \tensor \sfX + \nu_2 E^{2,1} \wedge E^{1,1} \tensor \sfX + \dots + \nu_8 \sum_{a=1}^m E^{1,a} \wedge E^{1,1} \tensor e_a{}^1$, and impose $0 = \sfX \cdot \kappa_4$ using the actions given in Section~\ref{S:FM}. Namely, $\sfX \cdot \sfY = \sfH$, $\sfX \cdot \sfH = -2\sfX$, and $\sfX \cdot e_a{}^b = 0$. Also, $\sfX \cdot E_{i,a} = E_{i+1,a}$, and so $\sfX \cdot E^{i,a} = -E^{i-1,a}$. We find that $0 = \sfX \cdot \kappa_4$ is equivalent to
\begin{align*}
& \nu_3 = \nu_4 = -\frac{\nu_1 + \nu_2}{2}, \qquad
 \nu_5 = 0, \qquad
 \nu_6 = \nu_8 = -\nu_7.
\end{align*}
Setting $(\nu_1,\nu_2,\nu_8) = (\mu_1,\mu_2,\mu_3)$ then yields the result.\hfill $\qed$
\end{itemize}	 \renewcommand{\qed}{}
\end{proof}

\begin{Corollary}\label{cor:kappa4}
All parameters involved in an algebraic model $(\ff;\g, \fp)$ of ODE type from Proposition {\rm \ref{P:k}} for $\bbU = \bbA_2^{\tr}$ are uniquely determined.
\end{Corollary}

\begin{proof}
Recall from Table \ref{Tab:EP} that $\bbU = \bbA_2^{\tr}$ arises for $n \ge 3$. By Propositions \ref{P:DR}\,(b) and~\ref{P:k}\,(c), any algebraic model $(\ff;\g, \fp)$ of ODE type with $\kappa_H = \Phi_\bbU \in \bbU$, for $\Phi_\bbU$ from Table~\ref{Tab:lwv}, and $\dim \ff = \fS_\bbU$ has
\begin{align}\label{E:fA2}
\begin{split}
&\ff = \sn\bigl\{E_{n,a},\dots, E_{2,a},\widehat{E}_{1,1},E_{1,b}, \widehat{E}_{0,1},E_{0,b}, \sfX \colon 1\le a \le m,\, 2 \le b \le m \bigr\} \oplus \fann (\Phi_\bbU),
\end{split}
\end{align}
 where $\fann (\Phi_\bbU)$ was given in Table \ref{Tab:fa}, and
\begin{align*}
\widehat{E}_{1,1} := E_{1,1} + (n-2)\zeta \sfZ_1 \in \ff, \qquad \widehat{E}_{0,1} := E_{0,1} + \zeta \sfY \in \ff
\end{align*}
 for some $\zeta \in \R$. Curvature is $\kappa =\Phi_{\bbA_2^{\tr}} + \kappa_4$, for $\kappa_4$ given in \eqref{E:kappa4}, and $[\cdot,\cdot]_\ff = [\cdot, \cdot] - \kappa(\cdot,\cdot)$.

Let us now impose the Jacobi identity. We define
\begin{align*} 
\Jac^\ff(x,y,z) := [x,[y,z]_{\ff}]_{\ff} -[[x,y]_{\ff}, z]_{\ff} -[y, [x, z]_{\ff}]_{\ff}, \qquad \forall x,y,z \in \ff.
\end{align*}
We calculate
\begin{align*}
\big[\widehat{E}_{1,1},[E_{0,2},E_{3,1}]_\ff\big]_\ff &= -(n-2)^2\left( 2\zeta +\frac{3(2m+3)-n(4m+3)}{mn(n+1)+6} \right)E_{2,2},\\
\big[E_{0,2},\big[\widehat{E}_{1,1},E_{3,1}\big]_\ff\big]_\ff &= -(n-2)^2\left(3\zeta +\frac{3(3m+5)-n(5m+3)}{mn(n+1)+6} \right)E_{2,2}, \\
\big[\big[\widehat{E}_{1,1},E_{0,2}\big]_\ff,E_{3,1}\big]_\ff &= -\frac{(n-1)(n-2)^2(mn-3)}{mn(n+1)+6}E_{2,2},
\end{align*}
so that
\begin{align}\label{E:zeta}
\Jac^\ff\bigl(\widehat{E}_{1,1}, E_{0,2},E_{3,1}\bigr) = 0 \quad \text{implies}\quad \zeta= \frac{\bigl(2n-n^2-3\bigr)m + 3n-9}{mn(n+1)+6}.
\end{align}
 Continuing in a similar manner, we find that
\begin{align}
& \Jac^\ff\bigl(\widehat{E}_{0,1}, E_{1,2},E_{3,1}\bigr) = 0 \quad \text{ implies}\quad \mu_1 = \frac{6(n-1)(n-2)(m+1)}{mn(n+1)+6},\label{E:mu}\\
& \Jac^\ff\bigl(\widehat{E}_{1,1}, E_{2,2},E_{2,1}\bigr) = 0 \quad \text{ implies}\quad \mu_2 = -\frac{6(n-1)(m+1)\bigl(m\bigl(n^3 + n^2-6n + 6\bigr) + 6\bigr)}{(mn(n+1)+6)^2}.\nonumber
\end{align}
Using $\zeta$ and $\mu_1$ above, we then have
\begin{align}\label{E:mu3}
\Jac^\ff(\widehat{E}_{0,1}, E_{2,2},E_{3,1}) = 0 \quad \text{implies} \quad \mu_3= 1-n.
\end{align}
As claimed, the parameters $\zeta$, $\mu_1$, $\mu_2$, $\mu_3$ are uniquely determined functions of $(n,m)$.

 We remark that the remaining Jacobi identities for $\ff$ are necessarily satisfied because the existence of a submaximally symmetric ODE model in the $\bbA_2^{\tr}$-branch (see Table \ref{Tab:CM}) guarantees the existence of a corresponding algebraic model of ODE type. (Necessarily, this is equivalent to the one found above.)
\end{proof}

\subsection{Conclusion} \label{S:ProvRes}

Let us now complete the proofs for Theorems \ref{T:I} and \ref{T:Main1}. Fix an irreducible C-class module $\bbU = \bbB_4, \bbA_2^{\tr}$, or $\bbA_2^{\tf}$ in the effective part $\bbE$, and recall the respective lowest weight vectors $\Phi_\bbU \in \bbU$ from Table \ref{Tab:lwv}. By Propositions \ref{P:DR} and \ref{P:k}, the classification of algebraic models $(\ff;\g,\fp)$ of ODE type with $0 \not \equiv \img(\kappa_H) \subset \bbU$ and $\dim \ff = \fS_{\bbU}$ is given in Table \ref{Tab:AM}. This completes step (i) of the classification strategy given in Remark \ref{R:Strategy}.

\begin{table}[h]\renewcommand{\arraystretch}{1.3}
\centering
$\begin{array}{|c|c|c|c|} \hline
n& \begin{array}{c}
\text{Irreducible C-class}\\ \text{module} \bbU\subset \bbE
\end{array} &\ff & \kappa \\ \hline\hline
2 & \bbB_4 & \fa^{\Phi_\bbU} & \begin{cases}
\Phi_\bbU, & \text{over}\ \C \\
\pm \Phi_\bbU, & \text{over}\ \R
\end{cases} \\\hline
\ge 3 & \bbA_2^{\tr} & \begin{array}{c}
\ff \ \text{in} \ \eqref{E:fA2} \ \text{with}\\
\zeta \ \text{in} \ \eqref{E:zeta}
\end{array} & \begin{array}{c}
\Phi_\bbU + \kappa_4,\ \text{with}\ \beta =1,\ \kappa_4 \ \text{in}\ \eqref{E:kappa4},\\
\text{and}\ \mu_1, \ \mu_2, \ \mu_3 \ \text{in} \ \eqref{E:mu} \ \text{and} \ \eqref{E:mu3}.
\end{array} \\ \hline
\ge 2 & \bbA_2^{\tf} & \fa^{\Phi_\bbU} & \Phi_\bbU\\ \hline
\end{array}$
\caption{Classification of algebraic models of ODE type with $0 \not \equiv \img(\kappa_H) \subset \bbU$ and $\dim \ff = \fS_\bbU$.}	
\label{Tab:AM}
\end{table}

 We now turn to step (ii) of Remark \ref{R:Strategy} and discuss how the ODE model classification in Table~\ref{Tab:CM} is deduced from the abstract classification in Table~\ref{Tab:AM}. Using fundamental invariants described in Section~\ref{S: C-class}, we confirm that these ODE lie in the claimed branches. In \cite[Table~10]{KT2021}, the point symmetries were given for all of these models with the exception of the second~$\bbB_4$ model. (See below for this case.) We confirm submaximal symmetry dimensions and deduce the associated algebraic models. (The latter is immediate by uniqueness in the $\bbA_2^{\tr}$, $\bbA_2^{\tf}$ cases, as well as the~$\bbB_4$ case over $\C$.)

To complete the proof of Theorem \ref{T:Main1}, we establish point-inequivalence over $\R$ of the following submaximally symmetric $\bbB_4$ models:
\begin{align}\label{E:3ODE}
u_3^a= \frac{3u_2^1 u_2^a}{2u^1_1} \qquad \text{or} \qquad u_3^a= \frac{3u_1^1u_2^1 u_2^a}{1 + (u^1_1)^2} \qquad \text{for}\quad 1 \leq a \leq m.
\end{align}
The point symmetry algebra $\cS$ of the former is given in \cite[Table 10]{KT2021}. On $J^0(\R,\R^m)$, the distributions $\ker\bigl({\rm d}u^1\bigr)$ and $\ker({\rm d}t)$ are each $\cS$-invariant, so these determine $\cS$-invariant foliations by level sets of $u^1$ and $t$, respectively. Total differentiation of the former implies that the level set $\bigl\{u_1^1=0\bigr\} \subset J^1(\R,\R^m)$ is $\cS$-invariant. Hence, the prolonged action of $\cS$ on $J^1(\R,\R^m)$ is not locally transitive.

In contrast, we now establish that the latter ODE in \eqref{E:3ODE} has symmetry algebra that acts locally transitively on $J^1(\R,\R^m)$. Its point symmetries are (for $1 \le a \le m$ and $2 \le b \le m$)
\begin{align*}
&\partial_t, \qquad \partial_{u^a}, \qquad t \partial_{u^b}, \qquad u^a\partial_{u^b}, \qquad \bigl(t^2 + \bigl(u^1\bigr)^2\bigr)\partial_{u^b}, \qquad u^1\partial_t -t\partial_{u^1}, \\
& t\partial_t + u^1\partial_{u^1} + 2 \sum_{b=2}^m u^b\partial_{u^b}, \qquad \bigl(t^2-\bigl(u^1\bigr)^2\bigr)\partial_t + 2t\sum_{a=1}^m u^a \partial_{u^a}, \\
& tu^1\partial_t + \frac{1}{2}\bigl(\bigl(u^1\bigr)^2-t^2\bigr)\partial_{u^1} + u^1\sum_{b=2}^m u^b\partial_{u^b}.
\end{align*}
In particular over $\R$, transitivity immediately follows from prolonging some of them to $J^1(\R,\R^m\!)$:
\begin{align*}
 \partial_t, \qquad
 \partial_{u^a}, \qquad
 t\partial_{u^b} + \partial_{u^b_1}, \qquad
 u^1\partial_t - t\partial_{u^1} - \bigl(1 + \bigl(u^1_1\bigr)^2\bigr) \partial_{u^1_1} - u^1_1 \sum_{b=2}^m u^b_1 \partial_{u^b_1}.
\end{align*}
Thus, the symmetry algebras of \eqref{E:3ODE} are point-inequivalent, and hence the ODEs are point-inequivalent.

 An alternate method is to establish that the symmetry algebras are {\em abstractly non-isomorphic}. Indeed, for $m \geq 2$, their semisimple parts are respectively $\fsl_2 \times \fsl_2 \times \fsl_{m-1}$ and $\mathfrak{so}_{1,3} \times \fsl_{m-1}$. (The $m=1$ case was remarked in \cite[p.~18]{AMS2021}.) However, this requires more details, while our argument given above is more direct. Moreover:

 \begin{Remark} Considering invariant foliations also gives the added bonus of suggesting a {\em complex} point-equivalence between the two ODE systems in \eqref{E:3ODE}. If we regard the latter ODE in \eqref{E:3ODE} over $\C$, then on $J^0(\C,\C^m)$, we find two invariant foliations by levels sets of $u^1 + {\rm i} t$ and $u^1 - {\rm i} t$ respectively. The invariant foliations discussed above in the first case now suggest considering the following (complex) point transformation
 \begin{align*}
 \bigl(\widetilde{t}, \widetilde{u}^1, \widetilde{u}^2,\dots , \widetilde{u}^m\bigr) = \bigl(u^1 + {\rm i} t, u^1 - {\rm i} t, u^2,\dots , u^m\bigr).
 \end{align*}
 We can straightforwardly verify that its prolongation pulls back the former ODE in \eqref{E:3ODE} (written in tilded variables) to the latter ODE in \eqref{E:3ODE}.
\end{Remark}

This completes the proof of Theorem \ref{T:Main1}. Following the remarks preceding Theorem \ref{T:I}, we have also proven the remaining Theorem \ref{T:I}(b) since for vector ODEs \eqref{ODE} of C-class of order~$n+1 \ge 3$, we have $\fS = \fM-2 = \fS_{\bbB_4} = \fS_{\bbA_2^{\tf}}$ only when $(n,m)= (2,2)$. This completes our proofs for Theorems \ref{T:I} and \ref{T:Main1}.

\appendix

\section[Harmonic curvature as the lowest degree component of curvature]{Harmonic curvature as the lowest degree component\\ of curvature} \label{S:kappa}

Fix $G$ and $P$ as in Section~\ref{S:FM} and recall from Section~\ref{S:CG} some basic notions of Cartan geometries $(\cG \to \cE, \omega)$ of type $(G, P)$ associated to ODEs \eqref{ODE}. We formulate Theorem \ref{T:lh} below stating that the harmonic curvature $\kappa_H$ can be identified with the {\em lowest degree} component (with respect to the grading element) of the curvature $\kappa$. (We note that this is used in the proof of Proposition~\ref{P:k}, which is essential in proving Theorems \ref{T:I} and \ref{T:Main1}.)

\begin{Definition}
Let $(\cG \to \cE, \omega)$ be a Cartan geometry of type $(G,P)$, let $\rho\colon G \to \GL(V)$ be a $G$-representation, and $\rho \circ\iota\colon P \to \GL(V)$ its restriction, where $\iota \colon P \hookrightarrow G$ is the canonical inclusion. A {\em tractor bundle} is an associated vector bundle $\cG \times_P V$ with respect to the $P$-representation $\rho \circ \iota$. Given the adjoint representation $\rho = \Ad\colon G \to \GL(\g)$, the tractor bundle $\cA \cE := \cG \times_P \g$ is called the {\em adjoint tractor bundle} (see \cite[Section~1.5.7]{CS2009} for further details).
\end{Definition}

Using the Cartan connection $\omega$, the tangent bundle $T\cE$ can be identified with the bundle $\cG \times_P (\g/\fp) $. Then, the $P$-invariant quotient map from $\g$ onto $\g /\fp$ gives rise to the natural projection $\Pi\colon \cA \cE \to T\cE$. Using this identification, we can regard the curvature as $\kappa \in \Omega^2(\cE, \cA \cE)$, i.e., $\cA\cE$-valued 2-form on $\cE$ \cite[Proposition~1.5.7]{CS2009}.

\begin{Definition} \label{D:RN}
Given a Cartan geometry $(\cG \to \cE, \omega)$ of type $(G,P)$ with curvature $\kappa \in \Omega^2(\cE, \cA \cE)$. Then
\begin{itemize}\itemsep=0pt
\item [(a)] $\omega$ is called {\em regular} if $\kappa \in \bigl(\Omega^2(\cE, \cA \cE)\bigr)^1$, i.e., $\kappa(T^i\cE, T^j\cE) \subset \cA^{i+j+1}\cE$, $\forall i$, $j < 0$.
\item [(b)] $\omega$ is called {\it normal} if $\partial^\ast \kappa = 0$.
\item [(c)] If $\omega$ is both regular and normal, then the {\em harmonic curvature} is $\kappa_H := \kappa \mod \img(\partial^*)$, which is a section of $\cG \times_P \frac{\ker \partial^*}{\img \partial^*}$.
\end{itemize}
\end{Definition}
Then, we have the following result.

\begin{Theorem} \label{T:lh}
Fix $G$ and $P$ as in Section~{\rm \ref{S:FM}}. Let $(\cG \to \cE, \omega)$ be a regular, normal Cartan geometry of type $(G,P)$ whose curvature $\kappa \in \left(\Omega^2(\cE,\cA \cE)\right)^\ell$ for some $\ell \ge 1$, i.e., $\kappa \left(T^i\cE, T^j\cE\right) \subset \cA^{i+j+\ell}\cE$ for all $i,j <0$. Then the induced section $\gr_\ell (\kappa) \in \gr_\ell \left(\Omega^2(\cE,\cA \cE)\right)$ coincides with the degree $\ell$ component of the harmonic curvature $\kappa_H$. Consequently, $\kappa_H \equiv 0$ implies $\kappa \equiv 0$.
\end{Theorem}

\begin{proof}
The statement was proved in \cite[Theorem~3.1.12]{CS2009} for parabolic geometries. The same proof works for our {\em non-parabolic} Cartan geometries associated to vector ODEs \eqref{ODE} of order~${\ge 3}$.
\end{proof}

\section[A necessary condition for coclosedness of Phi\_\{A\_2\^{}\{tr\}\}]{A necessary condition for coclosedness of $\boldsymbol{\Phi_{\bbA_2^{\tr}}}$}\label{S:coclosed}

 From Section~\ref{S:Phitr}, our strategy for computing a $\fsl(W)$-lowest weight vector $\Phi_{\bbA_2^{\tr}} \in \bbA_2^{\tr}$ involves imposing coclosedness, i.e., $\partial^* \Phi_{\bbA_2^{\tr}} = 0$, where $\partial^*$ was defined in Section~\ref{S:CG}. By adjointness of~$\partial$ and $\partial^*$ with respect to the inner product $\langle \cdot, \cdot \rangle$ on cochains induced from Definition \ref{D:metric}, we have
\begin{align}\label{coclosed}
\partial^\ast \Phi_{\bbA_2^{\tr}} = 0 \quad \iff \quad \big\langle \Phi_{\bbA_2^{\tr}},
\partial \psi \big\rangle = 0, \qquad \forall \psi \in \g^*_- \tensor \g.
\end{align}
 In order to pin down $\Phi_{\bbA_2^{\tr}}$ in Proposition \ref{P:ne}, only a small part of the conditions in \eqref{coclosed} will be in fact required. In this section, we identify a key condition (see Lemma \ref{L:norm}) that is essential to the proof of Proposition \ref{P:ne}.

 Recalling
$\fann(\Phi_{\bbA_2^{\tr}})$ given in Table \ref{Tab:fa}, let us restrict attention to $\psi$ lying in the subspace below.

\begin{Lemma} \label{L:1-cochain}
Suppose that $\psi \in \fg_-^* \otimes \fg$ has bi-grade $(1,1)$, with $\sfX \cdot \psi = 0$ and $\fann(\Phi_{\bbA_2^{\tr}}) \cdot \psi = 0$. Then $\psi$ is a multiple of
\begin{align} \label{E:Psi}
\Psi:= -2E^{2,1} \tensor \sfX + E^{1,1} \tensor \sfH + E^{0,1} \tensor \sfY.
\end{align}
\end{Lemma}

\begin{proof}
 Any $\psi \in \fg_-^* \otimes \fg$ with bi-grade $(1,1)$ lies in the span of
\begin{align*}
E^{2,a} \tensor \sfX, \qquad E^{1,a} \tensor \sfH,\qquad E^{0,a} \tensor \sfY, \qquad E^{1,a} \tensor e_b{}^c, \qquad 1 \le a,b,c \le m.
\end{align*}
Since $\sfX \cdot E_{i,a} = E_{i+1,a}$, then $\sfX \cdot E^{i,a} = -E^{i-1,a}$. Imposing $\sfX \cdot \psi = 0$ forces $\psi$ to lie in the span of
\begin{align} \label{E:1fmX}
-2 E^{2,a} \tensor \sfX + E^{1,a} \tensor \sfH + E^{0,a} \tensor \sfY, \qquad 1 \le a \le m.
\end{align}
	
Let us now impose $\fann(\Phi_{\bbA_2^{\tr}}) \cdot \psi = 0$. Recall from Table \ref{Tab:fa} that $q_d = e_d{}^d-e_{d+1}{}^{d+1} + \delta_1{}^d\sfZ_2 \in \fann(\Phi_{\bbA_2^{\tr}})$ for $1 \le d \le m-1$. Let $\fh \subset \fsl_m$ denote the standard Cartan subalgebra consisting of diagonal trace-free matrices, and $\epsilon_a \in \fh^*$ the standard weights for $\fh$. Since $\sfZ_2 \cdot \psi = \psi$, then
\begin{align*}
q_d \cdot \psi = 0 \quad \text{for}\quad 1 \le d \leq m-1 \quad\iff\quad \psi \quad \text{has weight} -\epsilon_1.
\end{align*}
Since each element of \eqref{E:1fmX} has weight $-\epsilon_a$, then being of weight $-\epsilon_1$ implies that $\psi$ is a multiple of \eqref{E:Psi}.	
We note that \eqref{E:Psi} is annihilated by all off-diagonal elements $e_f{}^d \in \fann\bigl(\Phi_{\bbA_2^{\tr}}\bigr)$, since we have $f \ge 2$ and $e_f{}^d$ commutes with $\{\sfX, \sfH, \sfY \}$. This completes the proof.
\end{proof}

 In terms of $\Phi^{i,j}= \sum_{a=1}^m E^{i,1} \wedge E^{j,a} \tensor E_{i+j-1,a}$ defined in \eqref{E:Phi_ij}, and using $\partial$ \eqref{E:cohdef}, we get
\begin{align} \label{E:dPsi}
\begin{split}
\partial \Psi = -2 \sum_{k=0}^{n-1} \Phi^{2,k} + \sum_{k=0}^n (2k-n) \Phi^{1,k} + \sum_{k=1}^n k(n+1-k)\Phi^{0,k}.
\end{split}
\end{align}
We then have the following necessary condition, which will be used in the proof of Proposition~\ref{P:ne}.
\begin{Lemma}\label{L:norm}
Take $\Phi_{\bbA_2^{\tr}}$ defined in Proposition {\rm \ref{P:GS}}, i.e., $\Phi_{\bbA_2^{\tr}} = \sum_{i,j=0}^nc_{i,j}\Phi^{i,j}$ with $c_{i,j}$ satisfying~\eqref{E:GS}. Then
\begin{align*}
0 = \sum_{k=0}^{n-1} \frac{(n-k)(k+1)}{n(n-1)}(c_{k,2}-mc_{2,k})
+ \sum_{k=0}^{n} \frac{2k-n}{n} (mc_{1,k}-c_{k,1}) + \sum_{k=1}^n (mc_{0,k}-c_{k,0}).
\end{align*}
\end{Lemma}

\begin{proof}
We evaluate \eqref{coclosed} for $\psi = \Psi$ given in \eqref{E:Psi}. In preparation for this, note that from Definition \ref{D:metric}, we have $\langle E_{k,a}, E_{k,a} \rangle = \frac{k!}{(n-k)!}$ and $\langle E^{k,a}, E^{k,a} \rangle = \frac{(n-k)!}{k!}$, and so
\begin{align*}
 \big|\big|E^{2,1} \wedge E^{k,a} \tensor E_{k+1,a}\big|\big|^2 = \big|\big|E^{2,1}\big|\big|^2 \big|\big|E^{k,a}\big|\big|^2\big|\big|E_{k+1,a}\big|\big|^2 = \frac{(n-k)(k+1)(n-2)!}{2}.
\end{align*}
 Hence, by bilinearity of $\langle \cdot,\cdot \rangle$ and orthogonality of the basis elements for $\g$ (see Definition \ref{D:metric}), we have
\begin{align*}
\big\langle \Phi^{2,k}, \Phi_{\bbA_2^{\tr}} \big\rangle
&= \sum_{i,j=0}^n\sum_{a,b=1}^m c_{i,j}\big\langle E^{2,1} \wedge E^{k,b} \tensor E_{k+1,b}, E^{i,1} \wedge E^{j,a} \tensor E_{i+j-1,a} \big\rangle\\
&= \sum_{i,j=0}^n \sum_{a=1}^m \bigl(\delta_i{}^2\delta_j{}^k - \delta_i{}^k \delta_j{}^2 \delta_a{}^1\bigr) c_{i,j} \big|\big|E^{2,1} \wedge E^{k,a}\tensor E_{k+1,a}\big|\big|^2\\
&= \sum_{a=1}^m \bigl(c_{2,k} - c_{k,2} \delta_a{}^1\bigr) \frac{(n-k)(k+1)(n-2)!}{2} \\
&= (mc_{2,k} - c_{k,2} ) \frac{(n-k)(k+1)(n-2)!}{2}.
\end{align*}
Similarly, we have
\begin{align*}
\big\langle \Phi^{1,k}, \Phi_{\bbA_2^{\tr}} \big\rangle = (n-1)!(mc_{1,k}-c_{k,1}), \qquad
\big\langle \Phi^{0,k}, \Phi_{\bbA_2^{\tr}} \big\rangle = \frac{n!}{k(n+1-k)}(mc_{0,k}-c_{k,0}).
\end{align*}
 We use these relations and \eqref{E:dPsi} to evaluate $0 = \langle\partial\Psi,\Phi_{\bbA_2^{\tr}} \rangle$ and obtain the claimed result.
 \end{proof}

\subsection*{Acknowledgements}

The authors acknowledge the use of the DifferentialGeometry package in \textsc{Maple}. We also acknowledge helpful conversations with Boris Kruglikov, Andreu Llabres, and Eivind Schneider. The research leading to these results has received funding from the Norwegian Financial Mechanism 2014--2021 (project registration number 2019/34/H/ST1/00636), the Troms\o{} Research Foundation (project ``Pure Mathematics in Norway''), and the UiT Aurora project MASCOT, and this article/publication is based upon work from COST Action CaLISTA CA21109 supported by COST (European Cooperation in Science and Technology), \url{https://www.cost.eu}.

\pdfbookmark[1]{References}{ref}
\LastPageEnding

\end{document}